\documentclass[leqno,11pt]{amsart}
\usepackage{latexsym,amsfonts,amsmath}
\newtheorem{theorem}{Theorem}[section]
\newtheorem{cor}[theorem]{Corollary}

\newtheorem{lemma}[theorem]{Lemma}
\newtheorem{prop}[theorem]{Proposition}
\theoremstyle{definition}
\newtheorem{definition}{Definition}[section]
\theoremstyle{remark}
\newtheorem{remark}{Remark}[section]
\numberwithin{equation}{section}
\newcommand{\nc}{\newcommand}
\newcommand{\C}{{\mathbb C}} \newcommand{\R}{{\mathbb R}}
\newcommand{\Z}{{\mathbb Z}} 
 \newcommand{\Q}{{\mathbb Q}}
\newcommand{\CA}{\mathcal{A}}

\newcommand{\vol}{\operatorname{vol}}

\renewcommand{\c}{\mathfrak{c}}
\renewcommand{\a}{\mathfrak{a}}

\newcommand{\res}{\operatornamewithlimits{Res}}

\nc{\at}{\C[\ga]}
\nc{\rat}{\C_\A[\ga]}
\nc{\ga}{\mathfrak{a}}
\nc{\gac}{\ga_\C}
\nc{\gd}{\mathfrak{g}}

\nc{\s}{\mathfrak{s}}
\nc{\gv}{\mathfrak{v}}
\nc{\gt}{\mathfrak{t}}
\nc{\dga}{\mathfrak{a}^*}
\nc{\dgd}{\mathfrak{g}^*}
\nc{\dgt}{\mathfrak{t}^*}
\nc{\gma}{\Gamma_\ga}
\nc{\gmd}{\Gamma_\gd}
\nc{\gmt}{\Gamma_\gt}
\nc{\dgma}{\Gamma_\ga^*}
\nc{\dgmd}{\Gamma_\gd^*}
\nc{\dgmt}{\Gamma_\gt^*}
\nc{\scd}{{\sc scd}}
\nc{\gc}{\mathfrak{c}}
\nc{\A}{\mathfrak{A}}
\nc{\B}{\mathfrak{B}}
\nc{\spt}{\mathrm{Sec}}
\nc{\bA}{\mathrm{BInd}(\A)}
\nc{\bAF}{\mathrm{BInd}(\A, F)}
\nc{\bB}{\mathrm{BInd}(\B)}
\nc{\con}{\mathrm{Cone}}
\nc{\consing}{\mathrm{Cone}_\mathrm{sing}}
\nc{\gcp}{{\bar\gc}^\perp}
\nc{\comp}{U(\A)}
\nc{\tensor}{\otimes}
\nc{\ra}{\rightarrow}
\nc{\isom}{\cong}
\nc{\dvolga}{\,d\mu_\Gamma^\ga}
\nc{\sumn}{\sum_{i=1}^n}
\nc{\sumr}{\sum_{i=1}^r}
\nc{\lr}[2]{\langle#1,#2\rangle}
\nc{\fs}[1]{F^\sigma_{(#1)}}
\nc{\epu}{\epsilon_{\mathrm{u}}}
\nc{\epl}{\epsilon_{\mathrm{l}}}
\nc{\bss}{\mathrm{Basis}}
\nc{\permr}{\mathcal{S}_r}
\nc{\sign}{\mathrm{sign}}
\nc{\mpl}{\mathrm{MP}_\lambda}
\nc{\flaga}{\mathcal{FL}(\A)}

\nc{\flag}{\mathcal{FL}}

\nc{\pone}{\mathbb{P}^1}
\nc{\cone}{\mathrm{cone}}
\nc{\ar}{\rightarrow}
\nc{\bl}{{\boldsymbol{\lambda}}}
\nc{\bgam}{\boldsymbol{\gamma}}
\nc{\bn}{\mathbf{n}}
\nc{\solp}{\mathrm{sol}_{\A,\bl}(\bn)}
\nc{\bfb}{\mathbf{b}}
\nc{\sigm}{\bar\sigma}
\nc{\hc}{H_\mathrm{c}}
\nc{\hbm}{H^{\mathrm{BM}}}
\nc{\prbl}{p^\R_\bl}
\nc{\gC}{\mathfrak{C}}
\nc{\vb}{t}
\nc{\ve}{\boldsymbol{\epsilon}}
\nc{\bq}{\mathbf{q}}

\newcommand{\eqs}{\operatorname{Eq}}
\nc{\Eq}{\mathrm{Part}}

\nc{\chuze}{\mathrm{Choose}}
\nc{\FF}{\mathcal{FL}(\A,\xi)}

\nc{\FFF}{\mathcal{FL^+}(\A,\xi)}
\nc{\iml}{\mathrm{im}(L)}
\nc{\trop}{t}
\nc{\tS}{\widetilde{\mathrm{Eq}}}
\nc{\sol}{\mathrm{sol}}
\nc{\ts}{\widetilde{\mathrm{sol}}}
\nc{\vkf}{\boldsymbol{\kappa}^F}
\nc{\tD}{\widetilde{D}}
\nc{\vsigma}{{\boldsymbol{\sigma}}}
\nc{\brho}{{\boldsymbol{\rho}}}
\nc{\const}{\mathrm{const}}
\nc{\cp}{{\mathrm{comp}}}
\nc{\ZZ}{\widehat{Z}}
\nc{\ZZF}{\widehat{Z}^F}
\nc{\JK}{\mathrm{JK}}
\nc{\overl}{\overset{\bl}<}

\title{Toric reduction and a conjecture of Batyrev and Materov}
\author{Andr{\'a}s Szenes and Mich{\`e}le Vergne} 
\begin{document}
\maketitle

\setcounter{section}{-1}
\section{Introduction}
\label{sec:intro}

This paper grew out of our efforts to understand the Toric Residue
Mirror Conjecture formulated by Batyrev and Materov in \cite{BM}.
This conjecture has its origin in Physics and is based on a work by
Morrison and Plesser \cite{MP}.  According to the philosophy of {\em
  mirror symmetry}, to every manifold in a certain class one can
associate a dual manifold, the so-called {\em mirror}, so that the
intersection numbers of the moduli spaces of holomorphic curves in one
of these manifolds are related to integrals of certain special
differential forms on the other.  While at the moment this mirror
manifold is only partially understood, there is an explicit
construction due to Victor Batyrev \cite{Bat1}, in which the two
manifolds are toric varieties whose defining data are related by a
natural duality notion for polytopes.

Let us recall the setting of the conjecture of Batyrev and Materov. 
Let $\gt$ be a $d$-dimensional real vector space endowed with an
integral structure: a lattice $\Gamma_{\gt}\subset\gt$ of full rank.
We denote by $\Gamma_{\gt}^*$ the embedded dual lattice
$\{v\in\dgt;\;\lr v\gamma\in\Z\text{ for all } \gamma\in\gt\}.$

Consider two convex polytopes, $\Pi\subset\gt$ and
$\check\Pi\subset\dgt$, containing the origin in their respective
interiors, and related by the duality
\[ \check \Pi=\{v\in\dgt;\;
\lr vb\geq-1\text{ for all }b\in\Pi\}.\] To simplify the exposition in
this introduction, we assume that both polytopes are simplicial and
have integral vertices. Then the correspondence between convex and
toric geometry associates to this data a pair of $d$-dimensional
polarized toric varieties: $V(\Pi)$ and $V(\check \Pi)$. Under our
present assumptions the polarizing line bundles, which we denote by
$L_\Pi$ and $L_{\check \Pi}$, are the anticanonical bundles of the
respective varieties. In the paper, we work in a more general setting
which is described in detail in \S\ref{sec:prelim}.

In the framework of the Batyrev-Materov conjecture, mirror symmetry
has two ``sides'': $A$ and $B$.  The $B$-side is characterized by a
certain function associated to the variety $V(\Pi)$ as follows. Each
point $\gamma$ of $\Pi\cap \Gamma_{\gt}$ gives rise to a holomorphic
section of $L_{\Pi}$. In particular, denote by $s_0$ the section
corresponding to the origin and by $\{s_i;\;i=1,\dots,n\}$ the
sections corresponding to the set of vertices $\{\beta_i;\;
i=1,\dots,n\}$ of $\Pi$.

Then, for a generic value of the complex vector parameter
$z=(z_1,\dots,z_n)$, the function $F_z=s_0-\sum_{i=1}^n z_i s_{i}$ is
a holomorphic section of $L_\Pi$, and the equation $F_z=0$ defines a
family of Calabi-Yau hypersurfaces in $V(\Pi)$ as $z$ varies. For each
$v\in \gt^*$, let $F_{v,z}=\sum_{i=1}^n \lr v{\beta_i} z_i s_i$, and
consider the ideal $I(z)$ generated by the sections $F_z$ and
$\{F_{v,z},\,v\in\gt^*\}$ in the homogeneous coordinate ring of
$V(\Pi)$. Then the toric residue introduced by Cox \cite{Cox} defines a
functional $\mathrm{TorRes}_{I(z)}$ on the space of sections of
$L_\Pi^{d}$, which vanishes on the subspace $H^0(V(\Pi),L^{d}_\Pi)\cap
I(z)$. Every homogeneous polynomial $P$ of degree $d$ in $n$ variables
gives rise to a section $S(P,z)=P(z_1 s_{1},\ldots,z_n s_{n})\in
H^0(V(\Pi), L_\Pi^{d})$, and thus we obtain a function
\begin{equation}\label{pb_intro}
 \langle P\rangle _{\B}(z)=\mathrm{TorRes}_{I(z)} S(P,z),
\end{equation}
 which is known to depend rationally on $z$. 
 
 Now we turn to the $A$-side of mirror symmetry, which is
 characterized by the solution to an enumerative problem on the
 variety $V(\check\Pi)$.  Introduce the notation
 $\ga=H_2(V(\check\Pi),\R)$, $\Gamma_\ga=H_2(V(\check\Pi),\Z)$ and
 also $\dga=H^2(V(\check\Pi),\R)$,
 $\Gamma^*_\ga=H^2(V(\check\Pi),\Z)$.  Recall from the theory of toric
 varieties that to each vertex $\beta_i$ of $\Pi$, and thus to each
 facet of $\check\Pi$, one can associate an integral element
 $\alpha_i$ of the second cohomology group $\Gamma^*_\ga$, which
 serves as the Poincar\'e dual of a particular torus-invariant divisor
 in $V(\check\Pi)$. The class $\kappa=\sumn\alpha_i$ is the Chern
 class of the anticanonical bundle of the variety; it plays an
 important role in the subject.

 Let $\ga_\mathrm{eff}\subset\ga$ be the cone of effective curves.
 For each $\lambda\in \Gamma_{\ga}\cap \ga_\mathrm{eff}$, Morrison and
 Plesser introduced a simplicial toric variety $\mpl$, which is a
 compactification of the space of holomorphic maps
\[ \{\iota:\pone\rightarrow V(\check\Pi); \; \iota_*(\phi)=\lambda\},\]
where $\phi$ is the fundamental class of $\pone$ in $H_2(\pone,\Z)$.
Their construction also produces a top-degree cohomology class
$\Phi^P_\lambda$ of $\mpl$, whose construction is similar to that of
the class $P(\alpha_1,\alpha_2,\ldots,\alpha_n)$ above. Then we can form the
generating series
\[ \langle P\rangle _{\A}(z)= \sum_{\lambda\in
  \Gamma_\ga\cap \ga_\mathrm{eff}} \int_{\mpl}\Phi^P_\lambda
\prod_{i=1}^n z_i^{\lr{\alpha_i}\lambda}. \] 

The {\em toric residue mirror conjecture} of Batyrev-Materov states that
this generating series is an expansion of the rational function
$\langle P\rangle _{\B}(z)$ in a certain domain of values of the
parameter $z$.  The precise statement of the conjecture in our general
framework is given in Theorem \ref{main1} after the preparations of
\S\S\ref{sec:prelim},\ref{sec:intersec},\ref{sec:MPmod}.
  
The main goal of the present paper is the proof of this theorem, however,
we feel that along the way we found a few results which are interesting
on their own right. Below we sketch these results, and, at the same
time, describe the structure of the paper and the highlights of the
proof.

After describing our setup and recalling the necessary facts from the
theory of toric varieties in \S\ref{sec:prelim}, we turn to the
intersection theory of toric varieties in \S\ref{sec:intersec}.
We approach the problem from the point of view of intersection numbers
on symplectic quotients initiated by Witten \cite{W} and Jeffrey and
Kirwan \cite{JK}. 

Let us consider an arbitrary simplicial toric variety $V$ of dimension
$d$. We maintain the notation we introduced for $V(\check\Pi)$:
$\alpha_i$, $i=1,\dots, n$, for the Poincar\'e duals of
torus-invariant divisors, and $\ga$, $\Gamma_\ga$, $\dga$,
$\Gamma^*_\ga$ for the appropriate second homology/cohomology groups.
We denote by $r$ the dimension of $\ga$. A polynomial $P$ of degree
$d$ in $n$ variables defines a top cohomology class of $V$ and one can
pose the problem of computing $\int_{V}P(\alpha_1,\dots,\alpha_n)$.
Witten \cite{W} and Jeffrey-Kirwan \cite{JK} gave rather complicated
analytic formulas for this quantity, involving some version of a
multidimensional inverse Laplace transform. In \cite{B-V} an algebraic
residue technique was given to compute these numbers; this algebraic
operation was named the Jeffrey-Kirwan residue.  Our first theorem,
Theorem \ref{A}, is a new iterated residue formula for the
Jeffrey-Kirwan residue, which maybe given the following homological
form.

Let $U=\{u\in \ga\tensor\C;\; \prod_{i=1}^n\alpha_i(u)\neq 0\}$ be the
complement of the complex hyperplane arrangement formed by the
zero-sets of the complexifications of the $\alpha$s in $\ga\tensor\C$,
and denote by $\gc$ the ample cone of $V$ in $\dga$.  For generic
$\xi\in\gc$ and a vector of auxiliary constants
$\ve=(\epsilon_1,\dots,\epsilon_r)$, we construct a cycle
$Z(\xi,\ve)\subset U(\A)$, which is a disjoint union
$\cup_{F\in\mathcal{FL}(\xi)} T_F(\ve)$ of oriented $r$-dimensional
real tori in $U$ indexed by a subset of flags of our hyperplane
arrangement depending on $\xi$.  Fix an appropriately normalized
holomorphic volume form $\dvolga$ on $\ga\tensor\C$. The integration
$f\mapsto\int_{T_F(\ve)} f\,\dvolga$ along one of these tori is called
an {\em iterated residue}; it is a simple algebraic functional on
holomorphic functions on $U(\A)$. Our integral formula (Theorem
\ref{A}) then takes the form
\[ \int_V P(\alpha_1,\dots,\alpha_n) = 
\int_{Z(\xi,\ve)}
\frac{P(\alpha_1,\dots,\alpha_n)\,\dvolga}{\alpha_1\dots\alpha_n}.
\]
Here, on the left hand side, we think of the $\alpha$s as cohomology
classes, while on the right hand side we consider them to be linear
functionals on $\ga\tensor\C$.

Next, in \S\ref{sec:MPmod}, we study the moduli spaces $\mpl$,
$\lambda\in\Gamma_{\ga}\cap\ga_\mathrm{eff}$, which are toric
varieties themselves. Using the results of \S\ref{sec:intersec}, we
derive an integral formula \eqref{sumpoly} for the generating function
$\langle P\rangle_\A(z)$ of the form
$\int_{Z(\xi,\ve)}P(u)\Lambda(u)$, where $\Lambda(u)$ is a meromorphic
top form in $U(\A)$. Here the constants $\ve$ need to be chosen
appropriately, in order to make sure $Z(\xi,\ve)$ avoids the poles of
$\Lambda$.

We turn to the $B$-side in \S\ref{sec:torres}.  We use a localized
formula \cite{CCD,CDS,BM} for the toric residue,
which has the form of a sum of the values of a certain rational
function over a finite set $O_\B(z)\subset V(\check\Pi)$. We make a
key observation (Proposition \ref{key} and Lemma \ref{cutout}) that
this finite set is naturally embedded into $U$ as the set of solutions
of the system of equations:
\[ \left\{ \prod_{i=1}^n \alpha_i(u)^{\lr{\alpha_i}\lambda} = \prod_{i=1}^n
z_i^{\lr{\alpha_i}\lambda}, \quad\lambda\in\Gamma_\ga\right\}.
\]
This infinite system of equations in the variable $u\in U$ easily
reduces to $r$ independent equations. This presentation of
$O_\B(z)$ allows us to write down an integral formula for $\langle
P\rangle_\B(z)$ in Proposition \ref{B}, which has the form of
$\int_{Z'} P(u)\Lambda(u)$, where $Z'$ is a another cycle in $U$
avoiding the poles of $\Lambda$. This way, we essentially reduce the
Batyrev-Materov conjecture to a topological problem of comparing
cycles.

The cycle $Z'$ is closely related to a real algebraic subvariety
$\ZZ(\xi) $of $U$
given by the set of equations
$$\ZZ(\xi)=\left\{u\in U(\A);\,\,
  \prod_{i=1}^n|\alpha_i(u)|^{\lr{\alpha_i}\lambda}=e^{-\langle \xi,
    \lambda\rangle } \text{ for all }\lambda\in \Gamma_\ga\right\}.$$
In \S\ref{sec:tropical} we prove the central result of the paper,
Theorem \ref{Zxi}, in which we compute the homology class of the cycle
$\ZZ(\xi)$ in $U$ for any generic $\xi$.  The proof  uses certain
type of degenerations reminiscent of the methods of tropical geometry
in real algebraic geometry (cf.  \cite{Viro,Sturm}).

In \S\ref{sec:completion} we specialize this result to the case when
the generic vector $\xi$ is near $\kappa=\sumn \alpha_i$, and
combining it with Theorem \ref{A}, we arrive at the statement that for
such $\xi$ the cycle $\ZZ(\xi)$ is contained in a small neighborhood
of the origin in $\ga\tensor\C$, and it is a small deformation of the
cycle $Z(\xi,\ve)$ which represents the Jeffrey-Kirwan residue (Theorem
\ref{C}). Armed with this result, the proof of the conjecture is
quickly completed.

We would like to end this introduction with a remark on the conditions
of our main result. Although here for simplicity we assumed that the
polytope $\Pi$ is simplicial and reflexive, neither of these
conditions are necessary. In the paper, we prove our result for an
arbitrary polytope with integral vertices, which contains the origin
in its interior.

Finally, we note that after this work was substantially completed, we
were informed by Lev A. Borisov that he had also obtained a proof of
the Toric Residue Mirror Conjecture by a completely different method.

\noindent {\sc Acknowledgments}. We would like to thank the
Mathematisches Forschungs\-institute Oberwolfach for wonderful working
conditions. We are thankful to Alicia Dickenstein for discussions on
toric residues, and to Victor Batyrev and Eugene Materov for
explaining their conjecture to us. The first author would like to
express his gratitude for the hospitality of Ecole Polytechnique and
acknowledge the support of OTKA.

\section{Preliminaries: Toric varieties}
\label{sec:prelim}

In this section, we describe standard facts from projective toric geometry.
The proofs will be mostly omitted  (cf. \cite{Danilov,Fulton,GKZ,HS}).

\subsection{Polytopes and toric varieties} \label{sec:poly_toric}
For a real vector space $\gv$ endowed with a lattice of full rank
$\Gamma_\gv$, denote by $\gv_\C$ the complexification
$\gv\otimes_\R\C$ of $\gv$, by $T_\gv$ the compact torus
$\gv/\Gamma_\gv$ and by $T_{\C\gv}$ the complexified torus
$\gv_\C/\Gamma_\gv$; finally, for $\gamma\in\Gamma^*_\gv$, where
\[
\Gamma_\gv^*=\{\gamma\in\gv^*;\; \lr\gamma v\in\Z\text{ for all
}v\in\Gamma_\gv\},\] denote by $e_\gamma$ the character $v\mapsto
e^{2\pi i\gamma(v)}$ of $T_\gv$. We will keep the notation $e_\gamma$
for the holomorphic extension of this character to the complexified
torus $\gv_\C/\Gamma_\gv$.  

Given a polytope $\Pi\subset\gv^*$ with integral vertices, one can
construct a polarized toric variety with action of the complex torus
$T_{\C\gv}$ as follows.  Consider the graded algebra
\begin{equation}
\label{gradedalg}
  \oplus \C e_\gamma g^k,k=1,\dots; \gamma\in
k\Pi\cap\Gamma^*_\gv,
\end{equation}
where the multiplication among the basis elements comes from addition
in $\gv^*$, and $g$ is an auxiliary variable marking the grading. This
algebra is the homogeneous ring of a polarized projective toric
variety $V(\Pi)$ endowed with a line bundle $L_\Pi$ and an action of
the torus $T_{\C\gv}$, which lift to an action on $L_\Pi$.  Each lattice point
$\gamma\in\Pi\cap\Gamma^*_\gv$ gives rise to a section $s_\gamma$ of
the line bundle $L_\Pi$, and the set $\{s_\gamma, \gamma\in
\Pi\cap\Gamma^*_\gv \}$ forms a basis of $H^0(V(\Pi),L_\Pi)$.  Note
that the toric variety $V(\Pi+t)$ corresponding to the polytope $\Pi$
translated by an element $t\in \Gamma^*_\gv$ is isomorphic to the
variety $V(\Pi)$, and the line bundle $L_\Pi$ is equivariantly
isomorphic to $L_\Pi\otimes \C_t$, where $\C_t$ is the one-dimensional
representation of $T_{\C\gv}$ corresponding to the character $e_t$.
Thus, starting from an integral  polytope in an affine space endowed
with a lattice, one can construct a well-defined polarized toric
variety.

\subsection{The quotient construction}
Now we give a more concrete description of toric varieties. Let
$\gd=\oplus_{i=1}^n\R\omega_i$ be an $n$-dimensional real vector space
with a fixed ordered basis, and let
\begin{equation}
\label{es}
 0\ra\ga\rightarrow\gd\overset\pi\rightarrow\gt\ra0
\end{equation}
be an exact sequence of finite dimensional real vector spaces of
dimensions $r,n$ and $d$, respectively. We assume that the lattice
$\gmd=\oplus_{i=1}^n\Z\omega_i$ intersects $\ga$ in a lattice
$\gma$ of full rank, and we denote the image $\pi(\gmd)$  in $\gt$ by
$\gmt$. This means that the sequence restricted to the lattices is
also exact. In this case the dual sequence
\begin{equation}
\label{des}
 0\ra\dgt\rightarrow\dgd\overset\mu\rightarrow\dga\ra0
\end{equation}
restricted to the dual lattices  
$\dgmt$, \,$\dgmd$ and $\dgma$, respectively, is also exact.

Denoting the elements of the dual basis by $\omega^i$, $i=1,\dots, n$,
we have $\dgd=\oplus_{i=1}^n\R\omega^i$; in particular, $\dgmd=
\oplus_{i=1}^n\Z\omega^i$.  Now introduce the notation $\alpha_i$ for
the image vector $\mu(\omega^i)$ in $\dgma$, $i=1,\dots,n$, and
consider the sequence $\A:=[\alpha_1,\alpha_2,\ldots,\alpha_n]$. We
emphasize that some of the $\alpha$s may coincide. The order of the
elements of this sequence will be immaterial, however.

\begin{definition}\label{defproj}
  We call a sequence $\A$ in $\dga$ {\em projective} if it lies in an
  open half space of the vector space $\dga$.
\end{definition}

The relevance of this condition will be explained below. Note that
according to our assumptions, the elements of $\A$ generate
$\Gamma^*_\ga$ over $\Z$, and this will {\em always} be tacitly
assumed in this paper.

\begin{definition} 
Let $\A$ be a not necessarily projective sequence
  in $\Gamma^*_\ga$.  Denote by $\bA$ the set of {\em basis index
    sets}, that is the set of those index subsets
  $\sigma\subset\{1,\dots,n\}$ for which the set
  $\{\alpha_i;\;i\in\sigma\}$ is a basis of $\dga$.  We will also use
  the notation
\[
\bgam^\sigma=(\gamma_1^\sigma,\dots,\gamma^\sigma_r)\] for the basis
associated to $\sigma\in\bA$; here a certain ordering of the basis
elements, say the one induced by the natural ordering of $\sigma$, has
been fixed. \end{definition}
\begin{definition}
For any set or sequence $S$ of vectors in a real vector space, denote by
$\con(S)$ the closed cone spanned by the elements of $S$.  Let us
consider the case of a projective sequence $\A$ in $\dga$.
We denote by $\consing(\A)$ the union of the boundaries of the
simplicial cones $\con(\bgam^\sigma)$, $\sigma\in\bA$.  Elements of
$\consing(\A)$ will be called {\em singular}, the others, {\em
  regular}. A connected component of $\con(\A)\setminus \consing(\A)$
is called a {\em chamber}. Then for a chamber $\gc$, we can define
$\mathrm{BInd}(\A,\gc)$ to be the set of those $\sigma\in\bA$ for
which $\con(\bgam^\sigma)\supset\gc$.\qed
\end{definition}

Now we assume that $\A$ is projective. Then we can proceed to
construct the toric variety $V_\A(\gc)$ as a quotient of the open set
\[U_\gc=\bigcup_{\sigma\in\mathrm{BInd}(\A,\gc)}
\left\{(z_1,\dots,z_n)\left|\, \prod_{i\in\sigma}z_i
    \neq0\right.\right\}\subset\C^n\] by the action of the
complexified torus $T_{\C\ga}$, where we let $T_{\C\ga}$ act on $\C^n$
diagonally with weights $(\alpha_1,\dots,\alpha_n)$. If $\A$ is
    projective then the quotient
\begin{equation}
\label{VAC}
V_\A(\gc)=U_\gc/T_{\C\ga}
\end{equation}
is a compact orbifold of dimension $d$.

To compare this construction to the one in\S\ref{sec:poly_toric}, let
$\theta$ be an integral point in $\gc$, and assume that the {\em
  partition polytope}
 \[
\Pi_\theta=\mu^{-1}(\theta)\bigcap\sum_{i=1}^n \R^{\geq0}\omega^i,
\]
which lies in an affine subspace of $\gd^*$ parallel to $\gt^*$, has
integral vertices. The toric variety $V_\A(\gc)$ is isomorphic to the
variety $V(\Pi_\theta)$ described before by its homogeneous ring. The
polarizing line bundle may be defined as $L_\theta=
U_\gc\times_{T_{\C\ga}}\C_\theta$, where $\C_\theta$ is the
one-dimensional representation of $T_{\C\ga}$ corresponding to the
character $e_\theta$. If
\[ \gamma=\sum_{i=1}^n \gamma_i \omega^i\in\Pi_\theta\cap \dgmd,\text{
  with }\gamma_i\in\Z^{\geq0},\,i=1,\dots,n, \]  then the
holomorphic function $\tilde s_\gamma:U_\gc\rightarrow  \C$ given
by $\tilde s_\gamma(z)=\prod_{i=1}^n z_i ^{\gamma_i}$ descends to a
section $s_\gamma$ of the line bundle $L_\theta$.

Note that under this quotient construction, for every $\eta\in\dgma$,
even for those not necessarily in $\gc$, a line bundle $L_\eta$ may be
defined by $L_\eta= U_\gc\times_{T_{\ga}}\C_\eta$. (In fact, $L_\eta$ is
usually only an orbi-bundle). The chamber $\gc$ is called the {\em
  ample cone} of the variety $V_\A(\gc)$ as the line bundles
corresponding to lattice points in $\gc$ are ample.

\subsection{Gale duality}
Let $\B$ be the sequence of vectors $\beta_i=\pi(\omega_i)\in\Gamma_\gt$,
where $\pi$ is the map in the exact sequence \eqref{es}. The sequence
$\B$ is called the {\em Gale dual} sequence to the sequence
$\A$. It immediately follows that taking the Gale dual of a sequence
twice, one recovers the original sequence.

The following Lemma describes the fundamental relation between
Gale dual vector configurations.

\begin{lemma}
\label{Gale}
A linear combination $\sum_{i=1}^n m_i\alpha_i$ vanishes if and only
if there is a linear functional $l\in\dgt$ such that $l(\beta_i)=m_i$.
\end{lemma}

This relation allows one to translate statements in the $\A$-language
into those in the Gale dual $\B$-language and vice versa. 

\begin{prop}
\label{collected}
  Let $\A$ be a projective sequence in $\dga$ and $\gc$ be
  chamber. Then
  \begin{enumerate}
  \item The Gale dual configuration $\B$ does not lie in any closed
    half space of $\gt$, that is $\sum_{i=1}^n\R^{\geq 0}\beta_i=\gt$.
\item  If $\sigma\in\bA$, then the complement
$\bar\sigma=\{1,\dots,n\}\setminus\sigma$
is an element of $\bB$.
\item Denote by $\bar{\bgam}$ the basis of $\gt$ corresponding to
  $\bar\sigma\in\bB$. The set of cones
  $\con(\bar{\bgam}^{\bar\sigma})$, $\sigma\in\mathrm{BInd}(\A,\gc)$
  forms a simplicial conic decomposition \scd$(\gc)$ of $\gt$.
\item The simplicial conic decomposition associated to the partition
  polytope $\Pi(\theta)$ coincides with \scd$(\gc)$ for any $\theta\in\gc$.
\end{enumerate}
  \end{prop}
  \begin{remark}
    A simplicial conic decomposition is also called a complete simplicial {\em
    fan}.
\end{remark}  

Now we prove a quantitative version of statement (2) of Proposition
\ref{collected}.  

Endow the vector spaces $\gd,\gt,\ga$ with orientations compatible
with the sequence \eqref{es}. Observe that a vector space $\gv$
endowed with a lattice of full rank and an {\em orientation} has a
natural translation-invariant volume form, that is an element of
$\wedge^{\dim\gv}\gv^*$, such that the signed volume of a unit
parallelepiped of the lattice is $\pm1$.  Accordingly, we have a
volume form on each vector space $\gd,\gt,\ga$; denote the volume form
on $\ga$ by $\dvolga$. Next,  for $\sigma\in\bA$, introduce the notation
$\vol_{\dga}(\sigma)$ for the signed volume of the parallelepiped
$\sum_{i\in\sigma} [0,1]\gamma_i^\sigma$. This means that we take the
volume of the parallelepiped measured in the units of the volume of a
basic parallelepiped of $\dgma$, and set the sign to $+1$ if the basis
$\bgam^\sigma$ is positively oriented, and to $-1$ otherwise;
$\vol_{\gt}(\bar\sigma)$ is defined similarly.  

Now we can formulate our first duality statement.
\begin{lemma}\label{volume} For $\sigma\in\bA$ we have
  $\bar\sigma\in\bB$, and
  $|\vol_{\dga}(\sigma)|=|\vol_{\gt}(\bar\sigma)|$.
\end{lemma}
\begin{proof}
  The exact sequence \eqref{es} gives rise to an isomorphism
$$I:\Lambda^d\gt\mapsto \Lambda^n\gd\otimes \Lambda^{r}\ga^*$$
as follows. For $y_1,y_2,\ldots,y_d\in \gt$ with representatives
$Y_1,Y_2,\ldots, Y_d$ in $\gd$, and $u_1,u_2,\ldots, u_r\in
\ga$, let
$$\langle  I(y_1\wedge y_2\wedge \cdots\wedge y_d),u_1\wedge u_2\wedge
\cdots \wedge u_r\rangle= Y_1\wedge Y_2\wedge \cdots \wedge
Y_d\wedge u_1\wedge u_2\wedge \cdots \wedge u_r.$$ Denote by
$\wedge\bgam^\sigma$ the form $\gamma_1^\sigma\wedge
\gamma_2^\sigma\wedge \cdots \wedge \gamma_r^\sigma$. Define
$\wedge\bar\bgam^{\bar\sigma}$ similarly. Then it is easy to
verify that $I(\wedge\bar\bgam^{\bar\sigma})=\pm (\omega_1\wedge
\omega_2\wedge \cdots \wedge \omega_n) \otimes \wedge
\bgam^\sigma.$ This implies the statement of the lemma.
\end{proof}

Now we translate a few important properties of vector configurations
into Gale dual language.

\begin{definition}
  Given a projective sequence $\A=[\alpha_i]_{i=1}^n$ in
  $\Gamma^*_\ga$, introduce the notation $\kappa=\sumn\alpha_i$. We
  call the sequence $\A$ {\em spanning} if for every
  $k\in\{1,\dots,n\}$ the vector $\kappa$ may be written as a
  non-negative linear combination $\kappa=\sum_{i=1}^nt_i\alpha_i$
  with $t_k=0$, and $t_i\neq0$, for $i=1,2,\dots,k-1,k+1,\dots,n$.
\end{definition}

\begin{lemma} \label{convexb} A sequence $\B$ is the set of vertices
  of a convex polytope containing the origin in its interior if and
  only if the Gale dual sequence $\A$ is projective and spanning.
\end{lemma}
\begin{proof} Indeed, for $\B=\{\beta_1,\dots\beta_n\}$ to be the set
  of vertices of a convex polytope is equivalent to the existence of
  linear functionals $h_k\in\dgt$ for $k=1,\dots,n$, such that
  $\lr{h_k}{\beta_k}=-1$ and $\lr{h_k}{\beta_i}>-1$ for $i\neq
  k$. Then according to Lemma \ref{Gale}, we have 
\[\sumn (\lr{h_k}{\beta_i}+1)\alpha_i = \kappa, \]
and this is exactly the spanning property for $\A$.
\end{proof}
\begin{remark}\label{extends}
  It is easy to see that if $\A$ is spanning, then the property
  described for $\kappa$ extends to any $\theta$ which is in a chamber
  $\gc$ containing $\kappa$ in its closure, i.e.  for every such
  $\theta$ and for each $k\in\{1,\dots,n\}$ one can find a
  non-negative integral linear combination
  $\theta=\sum_{i=1}^nt_i\alpha_i$ with $t_k=0$.
\end{remark}

Now we formulate two consequences of the spanning property. Recall
that according to statement (3) of Proposition \ref{collected}, the
 set of one-dimensional faces of the fan \scd$(\gc)$ is a {\em
  subset} of the set of rays $\{\R^{\geq 0}\beta_i;\;i=1,\dots,n\}$.
Also, note that there is a natural map $\chi_\gc:\Gamma_\dga\ar
H^2(V_\A(\gc),\Q)$ which associates to each lattice point $\theta$ the
first Chern class of the orbi-line-bundle $L_\theta$.
\begin{prop}
  \label{spanning} Assume that $\A\subset\dga$ is a spanning, projective
   sequence, and let $\gc$ be a chamber which contains $\kappa$ in its
   closure. Then 
  \begin{enumerate}
  \item the set of one-dimensional faces of the fan \scd$(\gc)$ is the set
  of rays   $\{\R^{\geq 0}\beta_i;\;i=1,\dots,n\}$,
\item the characteristic map $\chi_\gc$ is an isomorphism over $\Q$.
  \end{enumerate}
\end{prop}

The first statement follows from statement (4) of Proposition
\ref{collected}. Indeed, according to the above remark, for any
$\theta\in\gc$, the partition polytope $\Pi_\theta$ has exactly $n$
facets. The $k$th facet, which corresponds to the linear combinations
mentioned in the definition, is perpendicular to the Gale dual vector
$\beta_k$.  In the non-spanning case such a facet may disappear.

Next we can describe Batyrev's mirror dual toric varieties, which have 
respective actions of the tori $T_{\C\gt}$ and $T_{\C\gt^*}$.

Let $\B\subset\Gamma_\gt$ be the set of vertices of a convex polytope
$\Pi^\B$ containing the origin in its interior.  Then, on the one
hand, we can use this polytope to construct a projective toric variety
$V(\Pi^\B)$. On the other hand, consider {\em star-like}
triangulations of $\Pi^\B$, that is triangulations $\tau$ of $\Pi^\B$
with vertices at $\B\cup\{0\}$ such that every simplex contains the
origin. Clearly, such a triangulation $\tau$ gives rise to a
simplicial fan whose cones are the cones of the simplices of $\tau$
based at the origin.
\begin{prop}\label{starlike}
  Let $\B$ be a sequence of vectors whose elements serve as the
  vertices of an integral polytope $\Pi^\B$. Then the fan \scd$(\gc)$ of
  a chamber $\gc$ of the Gale dual configuration $\A$ induces a 
  star-like triangulation of $\Pi^\B$ if and only if $\gc$ contains the vector
  $\kappa=\sum_{i=1}^n\alpha_i$ in its closure.
\end{prop}

To summarize: the polytope $\Pi^\B$ corresponds to a toric variety
$V(\Pi^\B)$ on the one hand. On the other, it gives rise to a family
of ``mirror dual'' toric varieties $V_\A(\gc)$ corresponding to those
chambers $\gc$ of the Gale dual sequence $\A$ which contain $\kappa$
in their closure; the sequence $\A$ is spanning.

Finally, we recall the following definitions from
\cite{BM}.
Let $\check\Pi^\B$ be the
{\em dual polytope} of $\Pi^\B$ defined by
\[ \check\Pi^\B=\{l\in\dgt;\; \lr lb\geq-1,\,b\in\Pi^\B\}.\]
\begin{lemma}
  The dual polytope $\check\Pi^\B$ is a translate of the partition
  polytope $\Pi_\kappa$ associated to the Gale dual configuration $\A$.
\end{lemma}
\begin{proof} The point $t=\sum_{i=1}^n\omega^i$ is such that
$\mu(t)=\kappa$. It is easy to see that  $y\in\gd^*$ belongs to
$\Pi_\kappa$ if and only if $y-t\in\check\Pi^\B\subset\dgt$.
\end{proof}
\begin{definition}
  The polytope $\Pi^\B$ is called {\em reflexive} if the dual polytope
  $\check\Pi^\B$ has integral vertices.
\end{definition}

Batyrev and Materov consider dual pairs of reflexive polytopes. This
has the advantage of putting the toric variety and its mirror dual on
the same footing. In this paper, we will consider a more general
framework: we assume that $\B$ is the set of vertices of a polytope
with the origin in its interior, but no condition on the dual polytope
is imposed.

\section{Intersection numbers of toric quotients and the
  Jeffrey-Kirwan residue}
\label{sec:intersec}

In this section, $\A$ is any projective sequence in
$\Gamma_{\ga}^*$. Recall that we have chosen an orientation of $\ga$,
and that this, together with the lattice $\Gamma_\ga$ induces a
volume form $\dvolga$ on $\ga$.

Pick a chamber $\gc\subset\con(\A)$ and consider the orbifold toric
variety $V_\A(\gc)$. Since $V_\A(\gc)$ is a quotient
$U_\gc/T_{\C\ga}$, by the Chern-Weil construction, every polynomial
$Q$ on $\ga$ gives rise to a characteristic class $\chi(Q)$ of
$V_\A(\gc)$. Thus we have a Chern-Weil map
$\chi:\mathrm{Sym}(\ga^*)\rightarrow H^*(V_\A(\gc),\C)$ from the
polynomials on $\ga$ to the cohomology of $V_\A(\gc)$.  In particular,
for $\eta\in \ga^*$ the Chern class of the orbi-line-bundle $L_\eta$
is $\chi(\eta)$.

It is natural to look for formulas for the intersection numbers
$\int_{V_\A(\gc)} \chi(Q)$, where, of course, only the degree $d$
component of $Q$ contributes.  To write down a formula, we recall the
notion of the Jeffrey-Kirwan residue in a form written down by Brion
and Vergne \cite{B-V}.  Define $\comp$ to be the complement in $\a_\C$
of the complex hyperplane arrangement determined by $\A$:
$$\comp=\{u\in \ga_\C ;\, \alpha(u) \neq 0 \text{ for all }\alpha \in
\A\},$$
where we extended the functionals $\alpha_i$ from $\ga$ to
$\ga_\C$.

\begin{remark}  1. Note that $\alpha(u)$ and $\lr\alpha u$ stand for the same
  thing; we use one form or the other depending on whether we consider
  $u$ a
  variable or a constant.  \\
  2. The constructions of this section depend on the set of elements
  of $\A$, and do not depend on the multiplicities. We are not going
  to reflect this in the notation, however.
\end{remark}

Denote by $\rat$ the linear space of rational functions on $\ga_\C$
whose denominators are products of powers of elements of $\A$. The
space $\rat$ is $\Z$-graded by degree; the functions in $\rat$ are
regular on $\comp$.

Of particular importance will be certain functions in $\rat$ of degree
$-r$: for every $\sigma\in\bA$ denote by $f_\sigma$ the fraction
$1/\prod_{i\in\sigma}\alpha_i$. We will call such fractions {\em
  basic}.  Every function in $\rat$ of degree $-r$ may be decomposed
into a sum of basic fractions and degenerate fractions; degenerate
fractions are those for which the linear forms in the denominator do
not span $\dga$.  Now having fixed a chamber $\gc$, we define a
functional $\JK_\gc$ on $\rat$ called the Jeffrey-Kirwan residue (or
{\em JK-residue}) as follows.  Let
\begin{equation}
  \label{jkdef}
 \JK_{\c}(f_\sigma) =
\begin{cases}
  {|\vol_{\dga}(\bgam^\sigma)|}^{-1},&\text{if }\gc\subset
  \con(\bgam^\sigma),\\ 0, &\text{if }\gc\cap
  \con(\bgam^\sigma)=\emptyset.
\end{cases}
  \end{equation}
  Also, set the value of the JK-residue of a degenerate fraction
  or that of a rational function of pure degree different from $-r$
  equal to zero.

The definition of the functional $\JK_{\gc}(\cdot)$ is vastly
over-determined, as there are many linear relations among the basic
fractions $f_\sigma$, $\sigma\in\bA$.

\begin{prop}[\cite{B-V}]
  \label{jk-wdef}
The definition in \eqref{jkdef} is consistent and defines a
functional $\JK_{\gc}$ on $\rat$.
\end{prop}

One can give a homological interpretation to the JK-residue as follows.

\begin{lemma}\label{hom_interp}
For each chamber $\gc$ there is an homology class $h(\gc)\in
H_r(\comp,\R)$ such that
\[ \JK_\gc(f) = \frac1{\left(2\pi\sqrt{-1}\right)^r}
\int_{h(\gc)}f\,\dvolga \text{ for every }f\in\rat, \]
where $\dvolga$ is the translation invariant holomorphic volume form
defined above.
\end{lemma}
\begin{proof}
The integral on the right hand side is well-defined since the form
$f\,\dvolga$ is closed. The statement follows from Poincar\'e duality
and the fact that $H^r(\comp,\R)$ is spanned by holomorphic
differential forms of the form $f_\sigma\,\dvolga$, $\sigma\in\bA$
(cf. \cite{OT, VS}). 
\end{proof}

The integration over $V_\CA(\c)$ may be written in terms of the
Jeffrey-Kirwan residue as follows.
\begin{prop}[\cite{B-V}]
  \label{JK}
  Let $\A$ be a projective sequence in $\Gamma^*_\ga$,
   $\gc$ be a chamber and
  $Q$ be a polynomial on $\ga$. Then we have
\begin{equation}
  \label{integralJK}
  \int_{V_{\A}(\gc)} \chi(Q)=
  \JK_{\c}\left(\frac{Q}{\prod_{i=1}^n\alpha_i}\right).
\end{equation}
\end{prop}

Combining this with Lemma \ref{hom_interp} we obtain the formula
\begin{equation}
  \label{integralJK1}
  \int_{V_{\A}(\gc)} \chi(Q)=
\frac1{\left(2\pi\sqrt{-1}\right)^r}
\int_{h(\gc)}\frac{Q\,\dvolga}{\prod_{i=1}^n\alpha_i}.
\end{equation}

The main result of this section, Theorem \ref{A} may be interpreted as
a natural construction of a smooth cycle in $\comp$, which represents
the class $h(\gc)$.

We start with a few important notations and definitions related to our
hyperplane arrangement.  \smallskip

Let $\flaga$ be the finite set of flags
\[F = [F_0=\{0\}\subset F_1\subset F_2\subset  \dots \subset F_{r-1}
\subset F_r=\ga^*],\;\dim F_j=j,\] such that $\A$ contains a basis of
$F_j$ for each $j=1,\dots,r$.  For each $F\in\flaga$, we choose, once
and for all, an ordered basis ${\bgam}^F=(\gamma_1^F,\dots,\gamma^F_r)$
of $\dga$ with the following properties:
\begin{enumerate}
\item $\gamma_j^F\in\Gamma^*_\ga\tensor\Q$, for $j=1,\dots,r$,
\item $\{\gamma^F_m\}_{m=1}^j$ is a basis of $F_j$ for $j=1,\dots,r$,
\item the basis $\bgam^F$ is positively oriented,
\item $d\gamma^F_1\wedge\dots \wedge d\gamma^F_r=d\mu_{\Gamma}^{\ga}$.
\end{enumerate}
To each flag $F\in \flaga$, one can associate  a linear functional
$\res_F$ on $\rat$, called an {\em iterated residue}.
 We consider the elements of the basis $\bgam^F$
as coordinates on $\ga$ and we use the simplified notation
$u_j=\gamma^F_j(u)$ for $u\in\ga_\C$.  Then any rational function
$\phi\in\rat$ on $\ga_\C$ may be written as a rational function
$\phi^F$ of these coordinates:
\[ \phi(u) = \phi^F(u_1,\dots,u_r).\]
We define the iterated residue associated to the flag $F\in\flaga$
as the functional $\res_F:\rat\ar\C$ given by the formula
$$\res_F\phi=
\res_{u_r=0}\,du_r\;\res_{u_{r-1}=0}\cdots
\res_{u_1=0}\,du_1\;\phi^F(u_1,u_2,\ldots,u_r),$$ where each
residue is taken assuming that the variables with higher indices
have a fixed, nonzero value. 

It is easy to see that this linear form on $\rat$ depends only on the
flag $F$ and the volume form $\,d\mu_\Gamma^{\ga}$, and not on the
particular choice of the ordered basis $\bgam^F$. In fact, this
operation has a homological interpretation which is given below.

Let $N$ be a positive real number. Denote by $U(F,N)\subset
\ga_\C$ the open subset of $\ga_\C$ defined by
\[
U(F,N)=\{u\in \ga_\C\,;\; 0<N|\gamma^F_j(u)|< |\gamma^F_{j+1}(u)|,
j=1,2,\ldots, r-1\}.
\]

 The following Lemma is straightforward and its proof will be omitted.
\begin{lemma}
  \label{ufn}
  There exist positive constants $N_0$ and $c_0$ such that for $N>N_0$
  we have
\begin{enumerate}
\item $U(F,N)\subset\comp$ for all $F\in\flaga$, and
\item the sets $U(F,N)$, $F\in\flaga$ are disjoint.
\item If $\alpha_i\in F_j$ and $\alpha_k\in F_{j+1}\setminus F_j$ for
  some $F\in\flaga$ and $j<r$, then for every $u\in U(F,N)$ the
  inequality   $|\alpha_k(u)/\alpha_i(u)|>c_0 N$ holds.
\end{enumerate}
\end{lemma}

 From now on, when using the constant $N$, we will
assume that $N>N_0$. Note that the set $U(F,N)$ depends on the
choice of the basis $\bgam^F$ made above, but we will not reflect
this dependence in the notation explicitly. When $F$ is fixed, we
use as before the simplified notation $u_j=\gamma^F_j(u)$ for
$u\in\ga_\C$ so we can write

\[
U(F,N)=\{u\in \ga_\C\,;\; 0<N|u_j|< |u_{j+1}|, j=1,2,\ldots,
r-1\}.
\]

Observe that the set $U(F,N)$ is diffeomorphic to $\R_{>0}^r\times
(S^1)^r$, thus the $r$th homology of $U(F,N)$ is $1$-dimensional.
Choose a sequence of real numbers
$\ve:0<\epsilon_1\ll\epsilon_2\ll\dots \ll\epsilon_r$, where
$\epsilon\ll\delta$ means $N\epsilon<\delta$. Define the
torus
\begin{equation}
  \label{zefepsilon}
T_F(\ve) = \{u\in\gac;\; |u_j|=\epsilon_j,\, j=1,\dots,r\}\subset
U(F,N)\subset \comp,
\end{equation}
oriented by the form $d\arg u_1\wedge\dots\wedge d\arg u_r$. It is
easy to see that this cycle is a representative of a generator of
the homology $H_r(U(F,N),\Z)$.

The homology class of this cycle in
$\comp$ depends only on the flag $F$ and not on the chosen
positively oriented basis $\bgam^F$ of $F$. 

\begin{definition}\label{def_hf}
  Denote the homology class of the cycle $T_F(\ve)$ in $H_r(\comp,\Z)$
  by $h(F)$.  This produces a map $h:\flaga\ar H_r(\comp,\Z)$.
\end{definition}

\begin{lemma}[\cite{asz_ir}]
  \label{VS}
  For any $\phi\in\rat$ we have
\[ \frac1{\left(2\pi\sqrt{-1}\right)^r}\int_{h(F)}   \phi\;\,d\mu_{\Gamma}^{\ga}
\,\,= \res_F \phi,\]
where by integration over $h(F)$ we mean integration over any cycle
representing it.
\end{lemma}

Our goal is to write the functional $\JK_\c$ as a signed sum of
iterated residues $\res_F$. This will allow us to write $\JK_\c(\phi)$
as an integral of $\phi\, d\mu_{\Gamma}^{\ga}$ over the union of
corresponding cycles. The flags entering our formula will
depend on the choice of an element $\xi\in\gc$. This element will
have to satisfy additional conditions of regularity that we formulate
below.

\begin{definition}\label{sumregular}
 Denote by $\Sigma\A$ the set of elements of $\ga^*$ obtained by
partial sums of elements of $\A$:
\[\Sigma\A=\left\{\sum_{i\in\eta}\alpha_i;\; \eta\subset\{1,\dots,n\}\right\}.
\]
For each subset $\brho\subset \Sigma\A$ which forms a basis of
$\dga$, write $\xi$ in this basis: $\xi=\sum_{\gamma\in \brho}
u^\brho_\gamma(\xi)\, \gamma$.  Then introduce the quantity
\[
\mathrm{min}^{\Sigma\A}(\xi)=\min\{|u^\brho_\gamma(\xi)|;\;
\brho\subset\Sigma\A, \,\brho\text{ basis of
}\dga,\,\gamma\in\brho\}.
\]
 An element $\xi\in \ga^*$ will be called {\em regular with
  respect to $\Sigma\A$} if $\min^{\Sigma\A}(\xi)>0$.  For
$\tau>0$ we say that $\xi$ is {\em $\tau$-regular with respect to
  $\Sigma\A$}, if $\min^{\Sigma\A}(\xi)> \tau$.
\end{definition}
One could also say that a vector $\xi \in \ga^*$ is regular with
respect to $\Sigma\A$ if $\xi$ does not belong to any hyperplane
generated by elements of $\Sigma\A$. Sometimes, we will use the term
{\em sum-regular} for such vectors. Clearly, sum-regular vectors form
a dense open subset in $\dga$.

 Each flag $F\in\flaga$ introduces a partition of the elements of
     the sequence $\A=(\alpha_1,\dots,\alpha_n)$ induced by the
     representation of the space $\dga$ as a disjoint union
     $\cup_{j=1}^r F_j\setminus F_{j-1}$. For $j=1,\dots,r$, introduce
     the vectors
\begin{equation}\label{kappaj}
\kappa^F_j=\sum\{\alpha_i;\;i=1,\dots,n,\, \alpha_i\in F_j\}.
\end{equation}

  Note that the vectors $\kappa_j^F$ are in $\Sigma\A$, and that
$\kappa^F_r=\kappa=\sumn\alpha_i$ independently from $F$.

\begin{definition} \label{goodflags}
1. A  flag $F$ in $\flaga$ will be called {\em proper} if the
elements $\kappa_j^F,j=1,\ldots,r$
  are linearly independent.\\
2.  For each $F\in \flaga$, define a number $\nu(F)\in\{0, \pm 1\}$ as
follows: 
\begin{itemize}
\item set $\nu(F)=0$ if $F$ is not a proper
  flag;
\item if $F$ is a proper flag, then $\nu(F)$ is equal to $1$ or $-1$
  depending on whether the ordered basis
  $(\kappa_1^F,\kappa_2^F,\ldots,\kappa_r^F)$ of $\ga^*$ is positively
  or negatively oriented.
\end{itemize}
3. For a proper flag $F\in \flaga$, introduce the closed simplicial
cone $\s^+(F,\A)$ generated by the non-negative linear combinations of
the elements $\{\kappa^F_j, j=1,\ldots,r\}$:
\[ \s^+(F,\A) = \sum_{j=1}^r \R^{\geq 0}\kappa^F_j.\]
Then for $\xi\in \ga^*$,  denote by $\FFF$ the set of flags $F$ such
  that $\xi$ belongs to the cone $\s^+(F,\A)$.\qed
\end{definition}

Observe that if $\xi$ is sum-regular, then every flag $F\in\FFF$ is
proper, and thus for such $F$ we have  $\nu(F)=\pm 1$.

Now we are ready to formulate the main result of this section.

\begin{theorem}\label{A} 
  Let $\gc$ be any chamber of the projective sequence $\A$, and let
  $\xi$ be a vector in $\c$ which is regular with respect to
  $\Sigma\A$.  Then for every $\phi\in \rat$
 \begin{equation}
   \label{jkit}
\JK_{\gc}(\phi)=\sum_{F\in \FFF} \nu(F)
 \res_{F}\phi.
 \end{equation}
\end{theorem}

\begin{proof}
  Let $\sigma\in\bA$, and consider the basic fraction
$$f_{\sigma}=\frac{1}{\prod_{j=1}^r\gamma^\sigma_j}.$$

First, observe that it is sufficient to prove the theorem for
these basic fractions: $\phi=f_\sigma$ for $\sigma\in\bA$. Indeed,
both the Jeffrey-Kirwan residue and the iterated residues are
degree $-r$ operations on $\rat$. This allows us to restrict
$\phi$ to be of degree $-r$.  Now it is easy to check that the
iterated residues vanish on degenerate fractions, i.e. on
fractions whose denominators do not contain linear forms spanning
$\dga$. The JK-residue vanishes on degenerate fractions by
definition.

Thus we can assume $\phi=f_\sigma$. By the definition of the
chambers, the condition $\c\subset \con(\bgam^\sigma)$ is
equivalent to the condition $\xi\in \con(\bgam^\sigma)$.  Then
according to the definition of the JK-residue given in
\eqref{jkdef}, we have
$$\JK_{\c}(f_\sigma)= \begin{cases}
    \displaystyle\frac{1}{|\vol_{\dga}(\bgam^\sigma)|}, & \text{
    if } \c\subset \con(\bgam^\sigma)\\ 0, & \text{ otherwise.}
\end{cases}.$$

Now we compute the right hand side of \eqref{jkit} for
$\phi=f_\sigma$.  It is not hard to see that $\res_{F}(f_\sigma)$
is equal to $0$ unless the flag $F$ is such that its
$j$-dimensional component $F_j$ is spanned by elements of
$\bgam^\sigma$.  In other words, we have $\res_{F}f_\sigma\neq0$ for
some $F\in\flaga$ if
and only if $F$ is of the form 
\[ F^\pi(
\sigma)=(F_1^\pi(\sigma),\dots,F_r^\pi(\sigma))\;\text{with}\;
F_j^\pi(\sigma)= \Sigma_{k=1}^j \C\gamma^\sigma_{\pi(j)},
\] where $\pi$ is an element of $\Sigma_r$, the group of
permutations of $r$ indices. We will simply write $F(\sigma)$ in the
case when $\pi$ is the identity permutation. 

One can easily compute the appropriate iterated residue:
$$\res_{F^\pi(\sigma)}f_\sigma=\frac{(-1)^{\pi}}{\vol_{\ga^*}
(\bgam^\sigma)},$$
where we denoted by $(-1)^\pi$ the value of the alternating character
of $\Sigma_r$ on $\pi$.

Given an closed cone $C$, denote by $\chi[C]$ its characteristic
function. Using the above remarks, we can rewrite \eqref{jkit} as follows:
\begin{equation}
  \label{baricentric}
\sum_{\pi\in \Sigma_r}(-1)^{\pi}
\nu(F^\pi(\sigma))\chi[\s^+(F^\pi(\sigma),\A)](\xi)=
\chi[\con(\bgam^\sigma)](\xi)
\end{equation}
for any vector $\xi$ regular with respect to $\Sigma\A$.

As we will explain below in detail, this equality simply reflects the
subdivision into cones of $\con(\bgam^\sigma)$ based on the rays $\R^{\geq
  0}\kappa_j^{F^\pi(\sigma)}$, for $\pi\in\Sigma_r$, and
$j=1,\dots,r$.  Denote by $I(\sigma,\xi)$ the expression on the left
hand side of \eqref{baricentric}. We prove that
$I(\sigma,\xi)=\chi[\con(\bgam^\sigma)](\xi)$ by induction on the
dimension of $\ga$.

Consider the $(r-1)$-dimensional space $F_{r-1}(\sigma)$, the sequence
$\A'=\A\cap F_{r-1}(\sigma)$, and the index set $\sigma'$ obtained
from $\sigma$ by omitting its largest element.  For $\pi\in
\Sigma_{r-1}$, again we denote by $F^\pi(\sigma')$ the flag associated
to the permuted basis.  To compute $I(\sigma,\xi)$, we first study
$I_r(\sigma,\xi)$, the sum over the set $\Sigma_{r-1}$ of permutations
of the first $(r-1)$ indices:
$$I_r(\sigma,\xi)= \sum_{\pi\in \Sigma_{r-1}}(-1)^{\pi}
\nu(F^\pi(\sigma))\chi[\s^+(F^\pi(\sigma),\A)](\xi).$$

Recall that for any $F\in\flaga$ we have $\kappa^F_r=\kappa
\overset{\mathrm{def}}=\sum_{i=1}^n \alpha_i$.  There are two
cases:

 \begin{enumerate}

 \item

 The element $\kappa\in F_{r-1}$.

 \item

 $\ga^*=F_{r-1}\oplus \R \kappa$.

 \end{enumerate}

 We define $\nu_r\in \{-1,0,1\}$ as follows.  In the first case
 $\nu_r=0$.  In the second case, we write $\nu_r=\pm 1$ depending on
the orientation of
 $(\gamma^F_1,\gamma^F_2,\ldots,\gamma^F_{r-1},\kappa)$.
 Then in case (1), the sum $I_r(\sigma,\xi)$ is  equal to $0$,
 while in case (2), the cone $\s^+(F^\pi(\sigma),\A)$ is 
 equal to $\s^+(F^\pi(\sigma),\A')+ \R^+ \kappa$.  Writing
 $\xi=\xi'+ t\kappa$, we have $$\sum_{\pi\in
   \Sigma_{r-1}}(-1)^{\pi}
 \nu(F^\pi(\sigma))\chi[\s^+(F^\pi(\sigma),\A)](\xi)= $$
 $$
\begin{cases}
  \nu_r \sum_{\pi\in \Sigma_{r-1}}(-1)^{\pi}
  \nu(F^\pi(\sigma'))\chi[\s^+(F^\pi(\sigma'),\A')](\xi')  ,
&\text{ if } s> 0,\\
  0, & \text{ if } s<0.
\end{cases}$$

As $\xi$ is sum-regular, we cannot have $s=0$. Thus if $s>0$, then
the point $\xi'$ is sum-regular  with respect to $\A'$, and by the
induction hypothesis we conclude that $I_r(\sigma,\xi)=
\nu_r\chi[\con (\sigma')](\xi')$; if $s<0$, then we have
$I_r(\sigma,\xi)=0$.

Consider the closed cone $\con(\bgam^{\sigma'}\cup\{\kappa\})$.  The
preceding relation reads as
$$I_r(\sigma,\xi)= \nu_r\chi[\con(\bgam^{\sigma'}\cup\{\kappa\})](\xi).$$

It remains to sum over all circular permutations. Taking care of
the signs of the circular permutation and of orientations,
  we obtain this formula for the full sum:
$$I(\sigma,\xi)=\sum_{i=1}^r(-1)^i\nu_i\chi[C_i](\xi)$$
where $C_i$ is the cone generated by
$\bgam^\sigma\setminus\{\gamma^\sigma_i\}$ and $\kappa$, and
$\nu_i=\pm1$ depending on the orientation of this basis:
$(\gamma^\sigma_1,\gamma^\sigma_2,\dots,
\gamma^\sigma_{i-1},\gamma^\sigma_{i+1}, \dots, \gamma^\sigma_r,
\kappa)$. The fact that this sum equals $\chi[\con(\bgam^\sigma)](\xi)$ is
a straightforward exercise. This completes the proof of our theorem.
\end{proof}

\begin{remark} Using the results of \cite{asz_ir,B-V}, one can obtain
  a formula for the Jeffrey-Kirwan residue via iterated residues, using
  the concept of diagonal bases introduced in \cite{asz_ir}. Our
  present formula is quite different; it is more symmetric and seems
  to be computationally more efficient as well.
\end{remark}

Theorem  \ref{A} has the following
\begin{cor}\label{integralbis} 1. The equality 
\[ h(\gc) = \sum_{F\in \FFF}\nu(F) h(F)\] 
holds in $H_r(\comp,\Z)$. \\
2. The class $h(\gc)\in H_r(\comp,\R)$ is integral; it may be
represented by a disjoint union of embedded oriented tori.
\end{cor}
The first statement is a homological rewriting of Theorem \ref{A},
while the second follows from the fact that $h(F)$ is represented by
the torus $T_F(\ve)$. 

Thus we can reformulate Theorem \ref{A} in a third, integral form as
follows. Define the cycle
\begin{equation} \label{defzxi}
 Z(\xi) = \cup_{F\in \FFF}\nu(F) T_F(\ve),
\end{equation}
where $\ve$ is a vector of appropriate positive constants. Then
$Z(\xi)$ is an embedded oriented submanifold of $\comp$ depending on a
set of auxiliary constants, and we have
\begin{equation} 
\JK_\gc(\phi) = \int_{Z(\xi)} \phi\,\dvolga
\end{equation}

\section{The Morrison-Plesser moduli spaces}
\label{sec:MPmod}

In this section we assume that $\A$ is projective and spanning, and the
chamber $\gc$ contains $\kappa$ in its closure.  

Then, according to Proposition \ref{spanning}, we have a natural isomorphism
$H_2(V_\A(\gc),\Q)\isom \Gamma_\ga\tensor_\Z\Q$. Introduce the {\em
  cone of effective curves}
\[  \gcp=\{\lambda\in\ga;\; \langle  \xi,\lambda\rangle  \geq 0,
\text{ for all } \xi\in\gc\}. \]
Following Morrison and Plesser, we associate to each integral point
$\lambda\in \gcp\cap\gma$ of the cone of effective curves a toric
variety $\mpl$ together with a cohomology class $\Phi_\lambda\in
H^*(\mpl,\Z)$, called the {\em Morrison-Plesser moduli space} and its
{\em fundamental class}. This is a variant of the space of holomorphic
maps of $\pone$ into a fixed Calabi-Yau subvariety of $V_\A(\gc)$ with
a fixed image $\lambda$ of the fundamental class of $\pone$ (cf.
\cite{BM}).

Assume first that $\lambda\in\gma$ is such that $\langle
\alpha_i,\lambda\rangle \geq0,\;i=1,\dots,n.$ Then the
Morrison-Plesser toric variety $\mpl$ is the toric variety represented
by the data $\A^\lambda$, consisting of repetitions of the linear
forms $\alpha_i$ in $\dga$: each $\alpha_i$ is repeated $\langle
\alpha_i,\lambda\rangle +1$ times. Thus the total number of elements
of $\A^\lambda$ is $\langle \kappa,\lambda\rangle +n$, and the
dimension of the resulting toric variety is $\langle
\kappa,\lambda\rangle +d$. The polarizing chamber is the same one,
$\gc$, as that of the original toric variety. Thus we have $\mpl =
V_{\A^\lambda}(\gc)$.  The fundamental class is given by
$\Phi_\lambda=\chi(\kappa^{\langle \kappa,\lambda\rangle })$, where we
used the notation of the previous section.  We are interested in
intersection numbers of the variety $\mpl$ of the following form.

Fix a polynomial $P$ of degree $d$ in $n$ variables, and think of
it as a function on $\gd$. Denote the restriction of $P$ to $\ga$
by $P|\ga$; effectively, this means substituting $\alpha_i$ for
the $i$th argument of $P$. Then having fixed a sequence $\A$ and a
chamber $\gc$ in $\dga$, define
\begin{equation}
  \label{adef}
\langle   P\rangle  _{\lambda,\A,\gc} =
\int_{\mpl}\Phi_\lambda\chi(P|\ga).
\end{equation}
Using (\ref{integralJK1}), we can write
\begin{equation} \label{in-jk}
 \langle   P\rangle  _{\lambda,\A,\gc} =
 \frac1{\left(2\pi\sqrt{-1}\right)^r}\int_{h(\gc)}\frac{P(\alpha_1,\dots,\alpha_n)\,
 \kappa^{\langle  \kappa,\lambda\rangle  }\;\,d\mu_\Gamma^{\ga}}
 {\prod_{i=1}^n\alpha_i^{\langle  \alpha_i,\lambda\rangle  +1}},
\end{equation}
where $h(\gc)$ is the homology class representing the JK-residue.

Pick a cycle $Z[\gc]$ representing the homology class
$h(\gc)$, which satisfies the condition
\begin{equation}
  \label{eq:kappa_lto}
Z[\gc]\subset\comp\cap\{u\in\ga_\C;\;|\kappa(u)|<1\}.
\end{equation}
The cycle $Z(\xi)$ introduced in \eqref{defzxi} will be suitable if
the auxiliary constants $\epsilon_1,\dots,\epsilon_r$ are chosen
sufficiently small.

Now note that the rational function under the integral sign in
\eqref{in-jk} has exactly the correct degree: $-r$. This implies that
we can replace $\kappa^{\langle \kappa,\lambda\rangle }$ in the
formula by $(1-\kappa)^{-1}$ as follows:
\begin{equation}
  \label{integral2}
  \langle   P\rangle  _{\lambda,\A,\gc} = \frac1{\left(2\pi\sqrt{-1}\right)^r}\int_{Z[\gc]}
\frac{P(\alpha_1,\dots,\alpha_n)\; \,d\mu_\Gamma^{\ga}}
{(1-\kappa)\prod_{i=1}^n\alpha_i^{\langle \alpha_i,\lambda\rangle
+1}},
\end{equation}
Indeed to compute the right hand side of
\eqref{integral2} on such cycle $Z[\gc]$, we can replace
$1/(1-\kappa)$ by its absolutely convergent expansion
$\sum_{l=0}^{\infty}\kappa^l$. Then only the power $\kappa^{\langle
  \kappa,\lambda\rangle }$ gives a nonzero contribution to the
integral.

Further, observe that the right hand side of \eqref{integral2} is
meaningful for any $\lambda\in\gma$. Thus we can use it as a {\em
  definition} of the left hand side even for the cases when the
condition $\langle \alpha,\lambda\rangle \geq0$ does not hold for all
$\alpha\in\A$.

\begin{definition}\label{plambda}
  For any $\lambda\in \gma$, define the rational function
  $p_\lambda\in\rat$ of homogeneous degree $\langle
  \kappa,\lambda\rangle$ by
  $$p_\lambda=\prod_{i=1}^n\alpha_i ^{\langle \alpha_i,\lambda\rangle
  }.$$
One can write this function as quotient of two
  polynomials: $p_\lambda=p^+_\lambda/p^-_\lambda$, where
\[ p^+_\lambda=\prod_{\lr{\alpha_i}\lambda>0}
\alpha_i^{\langle \alpha_i,\lambda\rangle }, \quad
p^-_\lambda=\prod_{\lr{\alpha_i}\lambda<0} \alpha_i^{-\langle
  \alpha_i,\lambda\rangle . }\]
\end{definition}
The functions $p_\lambda$ satisfy the
relation $p_{\lambda_1}p_{\lambda_2}=p_{\lambda_1+\lambda_2}$ for any
$\lambda_1,\lambda_2\in \gma$.

\begin{definition}\label{Plambda}
  For any $\lambda\in \gma$ and degree $d$ polynomial $P$, we define
\begin{equation} \label{in-jk2}
 \langle   P\rangle  _{\lambda,\A,\gc} =
 \frac1{\left(2\pi\sqrt{-1}\right)^r}\int_{Z[\gc]}
 \frac{P(\alpha_1,\dots,\alpha_n)\,\;\,d\mu_\Gamma^{\ga}}
 {(1-\kappa)p_{\lambda} \prod_{i=1}^n\alpha_i},
\end{equation}
where $Z[\gc]$ is any cycle satisfying (\ref{eq:kappa_lto}).
\end{definition}

In the case when the condition $\langle \alpha,\lambda\rangle \geq0$
does not hold for all $\alpha\in\A$, the numbers $\langle P\rangle
_{\lambda,\A,\gc}$ may be interpreted as intersection numbers on a
modified version of the pair $(\mpl,\Phi_\lambda)$. Our convention for
$ \langle P\rangle _{\lambda,\A,\gc}$ induces a definition of the
fundamental class $\Phi_\lambda$ in this general case, which coincides
with the one given by Morrison and Plesser.  For details cf.
\cite{BM,MP}.

The next observation is central for our computations.
\begin{prop}
  \label{notincone}
  For  $\lambda\in\gma\setminus\gcp$,
one has $\langle   P\rangle  _{\lambda,\A,\gc} = 0.$
\end{prop}
\begin{proof}
  We can assume $\langle \kappa,\lambda\rangle \geq 0$, since $
  \langle P\rangle _{\lambda,\A,\gc} $ vanishes for $\langle
  \kappa,\lambda\rangle < 0$ by degree considerations.  Then \eqref{integral2}
  may be rewritten as
\begin{equation}
  \label{integral3}
  \langle   P\rangle  _{\lambda,\A,\gc} =\JK_{\gc}\left(
\frac{P(\alpha_1,\dots,\alpha_n)p^-_\lambda \kappa^{\langle
\kappa,\lambda\rangle  }}
{p^+_\lambda\prod_{i=1}^n\alpha_i}\right)
\end{equation}
Observe that the expression in  \eqref{integral3} is a JK-residue
of a rational function, denote it by $\phi_\lambda$, whose poles
lie on the hyperplanes $\alpha_i=0$, {\em with} $\langle
\alpha_i,\lambda \rangle \geq 0$. Indeed, if $\lr{\alpha_i}\lambda
<0$, then $\alpha_i$ occurs in the denominator
$p_\lambda^+\prod_{i=1}^n\alpha_i $ with multiplicity $1$, thus it
is canceled by a factor in $p_\lambda^-$ in the numerator.

Now comparing \eqref{integral3} to the definition of the
Jeffrey-Kirwan residue in \eqref{jkdef}, we see that
$\JK_{\gc}(\phi_\lambda)\neq 0$ implies that $\c$ is contained in
the cone generated by those $\alpha_i$ which satisfy
$\lr{\alpha_i}\lambda\geq0$. Consequently, $\lambda$, as a linear
functional on $\dga$ is positive on $\gc$, which is exactly the
condition $\lambda\in\gcp$.\end{proof}

Next, following \cite{BM}, we write down a generating series of the
numbers $\langle P\rangle_{\lambda,\A,\gc}$ for $\lambda$ varying in
the dual cone $\gcp$.  To this end, introduce the notation $z^\lambda$
for the Laurent monomial $\prod_{i=1}^n z_i^{\langle
  \alpha_i,\lambda\rangle }$ for any element $\lambda\in\gma$ and
$z=\sum_{i=1}^n z_i\omega_i\in\gd$.  Note that the restriction of the
function $z^\lambda$ to $\ga$ is exactly the rational function
$p_\lambda=p^+_\lambda/p^-_\lambda$.

Then the generating series of intersection numbers in which we are
interested  has the form
\begin{equation}\label{generating}
\langle P\rangle _{\A,\gc}(z)=\sum_{\lambda\in\gma\cap\gcp}\langle
P\rangle  _{\lambda,\A,\gc}z^\lambda.
\end{equation}

The chamber $\gc$, and thus $\gcp$, might be quite complicated, but
using Proposition \ref{notincone}, we can rewrite the generating
function \eqref{generating} very simply.

First an auxiliary statement:
\begin{lemma}
\label{unimod} Let $C$ be a closed rational  polyhedral cone in a
half-space of a real vector space $\mathfrak{v}$ of dimension $r$
endowed with a lattice $\Gamma$ of full rank, and let $\kappa$ be
a nonzero vector in $C$. Then there exist vectors
$v_1,v_2,\dots,v_r$ in $\Gamma$ with the properties
\begin{itemize}
\item $\sum_{j=1}^r\Z v_i=\Gamma$,
\item $v_j\in C$, for $j=1,\dots,r$,
\item $\kappa\in\sum_{j=1}^r\R^{\geq0}v_j$.
\end{itemize}
\end{lemma}
\begin{proof}
Indeed, the cone $C$ has a decomposition into simplicial cones
generated by $\Z$-bases of $\Gamma$.
\end{proof}

Now we return to our setup.
\begin{definition}
  Given a chamber $\gc$, we will call a set of vectors
  $\{\lambda_1,\dots,\lambda_r\}\subset\gma$ a {\em $\gc$-positive basis} if
  the following conditions are satisfied:
\begin{itemize}
\item $\sum_{j=1}^{r} \Z\lambda_j = \gma$,
\item $\sum_{j=1}^{r} \R_{\geq0}\lambda_j \supset \gcp$,
\item $\langle  \kappa,\lambda_j\rangle  \geq0$ for $j=1,\dots,r$.
\end{itemize}
\end{definition}
Apply Lemma \ref{unimod} to the pair $\kappa\in\bar\gc$.
Then  taking the dual basis guarantees the existence of a $\gc$-positive
basis. We fix such a basis and denote it by $\bl$.

Now observe that according to Proposition \ref{notincone} and the
second property of a $\gc$-positive basis, we can replace the sum in the
definition \eqref{generating} by the sum over a simplicial cone:
\[
\langle   P\rangle  _{\A,\gc}(z)=\sum\langle P\rangle_{\lambda,\A,\gc}\,
z^\lambda,\quad\lambda\in\sum_{j=1}^r\Z^{\geq0}\lambda_j.
\]

\smallskip
\noindent{\sc Notation.} Assume that a basis $\bl$ has been fixed, and 
let $z\in\gd_\C$ such that $z_i\neq 0$ for $i=1,\dots,n$. Then we
introduce the simplified notation $p_j, \,p^\pm_j$ for
$p_{\lambda_j}$, $p_{\lambda_j}^\pm$, respectively, and denote
$z^{\lambda_j}$ by $q_j$.
\smallskip

Using the integral definition of \eqref{in-jk2}, we can write
\begin{multline} \label{geometric_ser}
\langle   P\rangle_{\A,\gc}(z)=
 \frac1{\left(2\pi\sqrt{-1}\right)^r}\sum\int_{Z[\gc]}\prod_{j=1}^r\frac{q_j^{l_j}}{p_j^{l_j}}\cdot
\frac{P(\alpha_1,\dots,\alpha_n)\;\,d\mu_\Gamma^{\ga}}
 {(1-\kappa) \prod_{i=1}^n\alpha_i},
\\ \text{where the sum runs over }l_j\in\Z^{\geq0},\,j=1,\dots, r.
\end{multline}
If the condition
\begin{equation}
\label{qjcond} |q_j|<\max_{u\in Z[\gc]}|p_j(u)|, \quad j=1,\dots,r
\end{equation}
is satisfied, then the series is absolutely convergent. 

In fact, together with Proposition \ref{notincone}, this integral
representation allows us to determine the domain of convergence of
$\langle P\rangle_{\A,\gc}(z)$ more precisely. 
For $\lambda\in\Gamma_\ga$, define 
\begin{equation}\label{defeps}
 \epsilon_\lambda=\max_{u\in Z[\gc]}|p_\lambda(u)|,
\end{equation}
and consider the set
\begin{equation}\label{defw}
 W(Z[\gc])=\{z\in(\C^*)^n;\; |z^\lambda|<\epsilon_\lambda \text{ for every
}\lambda\in \Gamma_\ga\cap\gcp\}.
\end{equation}
Since both $z^\lambda$ and $\epsilon_\lambda$ are multiplicative
in $\lambda$, the set $W[\gc]$ is already defined by a set of
inequalities of the form $|z^\lambda|<\epsilon_\lambda$, where
$\lambda$ runs through a finite subset of $\Gamma_\ga\cap\gcp$ which
generates it as a semigroup. In particular, $W(Z[\gc])$ is open.

\begin{lemma}
\label{convergence} 
The series $\langle P\rangle_{\A,\gc}(z)$ converges for all $z\in W[\gc]$.
\end{lemma}
\begin{proof}
  We can decompose $\gcp$ into simplicial cones generated by $\Z$-bases of
  $\Gamma_\ga$. The sum in each such cone will be a convergent
  geometric series thanks to the inequalities defining $W[\gc]$.
\end{proof}

Now return to the fact that if the conditions \eqref{qjcond} hold,
then the series (\ref{geometric_ser}) converges absolutely. As a
consequence, we can exchange the order of summation and integration in
(\ref{geometric_ser}). 
Then we can sum the resulting geometric series under the integral sign
and arrive at the following statement.
\begin{prop}
\label{aside} Let $Z[\gc]$ be a cycle in $\comp$ representing
$h(\gc)$ and satisfying (\ref{eq:kappa_lto}). 
Fix a $\gc$-positive basis
$\lambda_1,\lambda_2,\dots,\lambda_r\in\gma$, and assume that
$z\in(\C^*)^n$ is such that the inequalities $|q_j|<\max_{u\in
  Z[\gc]}|p_j(u)|$ hold, where $q_j=z^{\lambda_j}$. Then we have
\begin{equation}
  \label{sumpoly}
\langle   P\rangle  _{\A,\gc}(z) =
\frac1{\left(2\pi\sqrt{-1}\right)^r}\int_{Z[\gc]}
\frac{P(\alpha_1,\dots,\alpha_n)\,\prod_{j=1}^r p_j
\,\,d\mu_\Gamma^{\ga}} {(1-\kappa)
\prod_{i=1}^n\alpha_i\prod_{j=1}^r \left(p_{j}-{q_j}\right)}.
\end{equation}
\end{prop}

\section{An integral formula for toric residues}
\label{sec:torres}

We start with the data considered so far: the exact sequences
\eqref{es} and \eqref{des}, the resulting sequence $\A$, which we
assume to be projective and spanning, a chamber $\gc\subset\dga$
containing $\kappa$ in its closure, a polynomial $P$ in $n$
variables, and a point $z\in\gd$.  Using these,
we defined a series $\langle P\rangle _{\A,\gc}(z)$  in the previous
section, and analyzed its domain of convergence. This series is the
object that the Batyrev-Materov conjecture associates to the $A$-side
of mirror symmetry.

Now we look at the same data in the Gale dual picture. According to
Lemma \ref{convexb}, the Gale dual sequence $\B\subset\Gamma_\gt$
serves as the set of vertices of a convex polytope $\Pi^\B$ containing
the origin.  The object on the $B$-side of the Batyrev-Materov
conjecture is the rational function $\langle P\rangle_\B(z)$ defined
in \eqref{pb_intro} of the introduction which uses the {\em toric residue}
of Cox (cf.  \cite{Cox,BM}).  As suggested in \cite{BM}, rather than
applying the original definition, we will use a localized formula
for toric residues \cite{CCD,CDS} which we recall below in \eqref{pb}.
The applicability of this localization formula in our case was kindly
explained to us by Alicia Dickenstein.

Consider the function $f=1-\sumn z_ie_{\beta_i}$, parameterized by our
chosen point $z=\sumn z_i\omega_i\in\gd$. Pick a $\Z$-basis
$(h_1,\dots,h_d)$ of $\dgmt$ and form the toric partial derivatives
\[   f_k = -\sumn z_i\langle   h_k,\beta_i\rangle   e_{\beta_i},\quad
k=1,\dots,d,\]
which assemble into the toric gradient $\nabla f=(f_1,\dots,f_d)$.
We can go on and define the toric Hessian of the function $f$ as
\[ H_f = \det\left(\sumn
  \lr{h_j}{\beta_i}\lr{h_k}{\beta_i}z_ie_{\beta_i}\right)_{j,k=1}^d.\]

Now denote by $O_\B(z)$ the set of toric critical points of $f$,
i.e. the set
\[ O_\B(z) = \{\nabla f=0\}\subset T_{\C\dgt}. \]
For generic $z$, this set is discrete, and the toric critical points
are non-degenerate. We take the following version of toric residues
localized at the toric critical points to be the definition of
$\langle P\rangle _\B(z)$:
\begin{equation}\label{pb}
\langle   P\rangle  _\B(z)=
\sum \frac{\tilde P(w)}{f(w)H_f(w)},\quad w\in  O_{\B}(z),
\end{equation}
where the function $\tilde P$ is obtained by substituting
$z_ie_{\beta_i}$ for $x_i$ in our degree $d$ polynomial
$P(x_1,\dots,x_n)$. 

In our setup, the conjecture of Batyrev and Materov generalizes to the
following statement.
\begin{theorem}\label{main1}
  Let $\A$ be a projective, spanning sequence, and $\gc$ be a chamber
  whose closure contains $\kappa$. Choose a cycle $Z[\gc]$ in $\comp$
  representing $h(\gc)$ and satisfying (\ref{eq:kappa_lto}), and let $z\in
  W(Z[\gc])$, where $W(Z[\gc])$ is defined in (\ref{defw}).  Then the
  series $\langle P\rangle_{\A,\gc}(z)$ converges absolutely, moreover, we
  have
\begin{equation}
  \label{eq:main}
\langle   P\rangle_{\A,\gc}(z) = \langle   P\rangle  _{\B}(z).
\end{equation}
\end{theorem}

The proof of this theorem is given at the end of the paper in
\S\ref{sec:theproof}. Its main ingredients are Propositions
\ref{aside} and \ref{B}, and
Theorem \ref{C}, which, in turn, follows from Theorems \ref{A} and \ref{Zxi}.
\begin{remark}
  1. In the course of the proof, we will construct an explicit cycle
     $Z[\gc]$, thus the domain of convergence of
     $\langle P\rangle_{\A,\gc}(z)$ will also be given explicitly.\\
     2. We think of the right hand side here as a rational function
     of $z$ given by the toric residue.  Note that, in particular,
     the right hand side does not depend on the choice of the chamber
     $\gc$. This dependence is encoded in the
     domain of convergence. \\
     3. The conjecture in \cite{BM} is formulated for the case of
     toric varieties corresponding to reflexive polytopes. As
     explained at the end of \S\ref{sec:prelim}, this corresponds to
     the special case of the partition polytope $\Pi_\kappa$ having
     integral vertices.
\end{remark}

The key observation that begins relating the two seemingly unrelated
expressions in \ref{eq:main} is the following.  Take a point $z\in\gd$
with all its coordinates $z_i\neq 0$, and embed $T_{\C\dgt}$ into
$\gd_\C$ via the formula
\begin{equation}
  \label{embedding}
w\mapsto(z_1 e_{\beta_1}(w),\dots,z_ne_{\beta_n}(w)).
\end{equation}
This means that we consider the natural action of $T_{\C\dgt}$ on
$\gd_\C$ given by the set of weights $\B$, and look at the orbit
$\mathrm{Orb}_\B(z)$ of the point $z\in\gd_\C$. Note that the coordinates
of a point in $\mathrm{Orb}_\B(z)$ are also all nonzero.
\begin{prop}
  \label{key}
  Under the embedding \eqref{embedding}, the set of critical points
  $O_\B(z)$ corresponds to the intersection of the orbit
  $\mathrm{Orb}_\B(z)$ with the linear subspace $\ga_\C$ in $\gd_\C$.
\end{prop}
\begin{proof}: We need to show that if $\nabla f=0$, then $(z_1
e_{\beta_1},\dots,z_ne_{\beta_n})\in\ga_\C$. Since $\nabla f=0$
exactly when $\sumn z_i e_{\beta_i}\beta_i=0$, this immediately
follows from Lemma \ref{Gale}.\end{proof}

This Proposition makes contact between the dual toric variety and the
second homology of the original toric variety. What is more, clearly,
the functions that appear in the definition of $\langle
P\rangle_\B(z)$ in \eqref{pb} all come as restrictions of
functions from the ambient space $\gd_\C$ identified with $\C^n$ by
$(x_1,\dots,x_n)\mapsto\sum_{i=1}^n x_i \omega_i$. Indeed $P$ was a
polynomial in $n$ variables, $f$ is the restriction of the function
$1-\sumn x_i$ and the Hessian may be considered as the restriction of
the function
\begin{equation}
  \label{DBdef}
 D^\B(x) = \det\left(\sumn
  \lr{h_j}{\beta_i}\lr{h_k}{\beta_i}x_i\right)_{j,k=1}^d,
\end{equation}
which is a degree $d$ polynomial on $\gd$.

Now that we may think of $O_\B(z)$ as a finite subset of $\ga_\C$, we
would like to know something about its geometry.  Fix a $\Z$-basis
$\bl=(\lambda_1,\dots,\lambda_r)$ of $\gma$, not necessarily a
$\gc$-positive basis, and recall the notation
\[ p_j(u)=p_j^+(u)/p_j^-(u)=\prod_{i=1}^n
\alpha_i(u)^{\lr{\alpha_i}{\lambda_j}}, \;\text{ for }\;j=1,\dots,r.
\]
Also, having fixed an appropriate $z\in\gd_\C$, with $z_i\neq 0$ for
all $i$, again denote by $q_j$ the number $z^{\lambda_j}$. As
$\ga_{\C}$ is embedded in $\gd_\C$ by $u\mapsto ( \alpha_1(u)
,\ldots,\alpha_n(u) )$, the set $O_\B(z)$ is contained in $U(\A)$.
Thus our set of critical points of the function $f$ on $T_{\C\dgt}$
becomes a finite subset of $\comp$. As such, it is cut out from
$\comp$ by some equations.
\begin{lemma}
  \label{cutout}
We have
\[ O_\B(z) = \{ u\in U(\A); \; p_j(u)=q_j,\,j=1,\dots,r\}.
\]
\end{lemma}

\begin{proof}
  Considering Proposition \ref{key}, the statement follows if we show
  that the torus $T_{\C\dgt}$ embedded via \eqref{embedding} is cut
  out from $\gd_\C$ by the equations
\[ \prod_{i=1}^n x_i^{\lr{\alpha_i}{\lambda_j}}=q_j,\;j=1,\dots,r.
\]
Thus what we need to show is that if $\prod_{i=1}^n
x_i^{\lr{\alpha_i}{\lambda_j}}=1$, for $j=1,\dots,r$, then for some
$h\in\dgt_\C$ we have $x_i=e^{2\pi\sqrt{-1}\lr h{\beta_i}}$.
Representing $x_i$ as $e^{2\pi \sqrt{-1}l_i}$ and using the fact that
$\bl$ is a basis of $\Gamma_\ga$ over $\Z$, we see that $\sumn
l_i\alpha_i\in\dgma$. According to our assumptions, the $\alpha_i$s
generate the lattice $\dgma$ over $\Z$. Since the $l_i$s are defined
only up to integers, by choosing them appropriately, we may assume
that $\sumn l_i\alpha_i=0$. Then the basic property of Gale duality,
Lemma \ref{Gale}, completes the proof.
\end{proof}

We can summarize our results so far as follows: we have
\begin{equation}
  \label{summary}
 \langle   P\rangle  _\B(z) =
   \sum\frac{P(\alpha_1(u),\dots,\alpha_n(u))}
{(1-\kappa(u))D^\B(\alpha_1(u),\dots,\alpha_n(u))},
\end{equation}
where, as usual,  $\kappa=\sumn \alpha_i$, and the sum runs over the
finite set $\{ u\in U(\A); \; p_j(u)=q_j,\,j=1,\dots,r\}$.

The statement of Theorem \ref{main1} is thus reduced to showing that
the integral in \eqref{sumpoly} of a rational differential form, which
we denote by $\Lambda$, over the cycle $Z[\gc]\subset U( \A)\subset
\ga_\C$ is equal to the expression in \eqref{summary}: a finite sum of
the values of a rational function over a finite set of common zeros of
$r$ other rational functions.  The first step of the proof, completed
in this section, will be showing that this finite sum also has a
representation as an integral of the same form $\Lambda$ over a
different cycle. The second step, which will take up the
rest of the paper, will be showing the equivalence of the two cycles.

First, we compute the coefficients of the polynomial $D^\B(x)$ defined
in \eqref{DB} explicitly.
\begin{lemma}\label{DB}
We have
\begin{equation}
  \label{DBformula}
D^{\B}(x)=
\sum_{\bar\sigma\in\bB}\vol_{\gt}({\bar\bgam}^{\bar\sigma})^2
\prod_{i\in   \bar\sigma} x_i.
\end{equation}
\end{lemma}
\begin{proof}
  Thinking of the vectors $\beta_i$, as $d$-component column vectors
  written in the basis $\{h_k\}_{k=1}^d$, we can write the matrix
  $M(x)$ the determinant of which is $D^\B(x)$ as
\[ M(x) = \sumn x_i\beta_i\beta_i^T, \]
where $\beta^T$ is the transposed matrix: a row vector. Using the fact
that each of the terms in this sum is a rank-1 matrix, we can expand
$\det(M(x))$ as
\[  \det(M(x))= \sum_{\bar\sigma\in\bB}
  \det\left(\sum_{i\in\bar\sigma}x_i\beta_i\beta^T_i\right).
\]
The term of this sum corresponding to the basis $\bar\sigma\in\bB$,
written in the basis $\bar\sigma$ itself, is simply a diagonal matrix
with entries $\{x_i;\,i\in\bar\sigma\}$ on the diagonal. This
immediately implies \eqref{DBformula}. \end{proof}

Next we compute the Jacobian matrix of the vector valued function
\begin{equation}\label{defp}
 p=(p_1,\dots,p_r): \comp\ar\C^{*r}.
\end{equation}

\begin{prop}
  \label{DA}
  Define the rational function $D_\A$ on $\gd_\C$ as
\[ D_\A(x) =
\det\left(\sumn\frac{\lr{\alpha_i}{\lambda_l}
   \lr{\alpha_i}{\lambda_m}}{x_i}\right)_{l,m=1}^r.
\]
Then we have
\begin{enumerate}
\item $\displaystyle D_\A(x) = \sum_{\sigma\in\bA}\vol_{\dga}(\bgam^\sigma)^2
  \prod_{i\in\sigma}\frac1{x_i},$
 \item $ \displaystyle D_\A(x)\prod_{i=1}^nx_i=D^\B(x).$
\item $\displaystyle \frac{dp_1}{p_1}\wedge\dots\wedge \frac{dp_r}{p_r}(u) =
   D_\A(\alpha_1(u),\dots,\alpha_n(u))\; \,d\mu_\Gamma^{\ga}.$
\end{enumerate}
\end{prop}
\begin{proof}
  The proof of (1) is exactly the same as that of Lemma \ref{DB}. Then
  (1) and Lemma \ref{DB} together with Lemma \ref{volume} imply
  (2). Finally, (3) is a simple calculation: Taking the partial
  derivative of $p_l$ with respect to $\lambda_m$ is exactly $p_l$
  times the corresponding entry of the matrix in the definition of $D_\A(x)$.
\end{proof}
\begin{cor}
The map $p:\comp\ar \C^{*r}$ is generically nonsingular.
\end{cor}
Indeed, statements (1) and (3) of Proposition \ref{DA} compute the
Jacobian of this map explicitly.  Since $\A$ is projective, there is a
$u\in\comp$ such that $\alpha_i(u)>0$, $i=1,\dots,n$, and at such $u$
the sum in statement (1) is clearly positive.\qed

Now we are ready to present our residue formula for $\langle
P\rangle_\B(z)$. For $z\in(\C^*)^n$ and $\delta>0$ let
\[Z_\delta(\bl,q)=\{u\in U(\A);\; |p_j(u)-q_j|=\delta,\,j=1,\dots,r\}\]
oriented by the form $d\arg (p_1-q_1)\wedge\dots\wedge d\arg
(p_r-q_r)$.

\begin{prop}   \label{B} 
Let $z\in\C^{*r}$ be such that the set $O_\B(z)\subset\comp$ is finite and the
function $(1-\kappa)D^\B(\alpha_1,\dots,\alpha_n)$ does not vanish on
it, and let $U(z)$ be a small neighborhood of $O_\B(z)$ in $\comp$. 
Then for sufficiently small $\delta>0$, we have
\begin{equation}\label{pb_integral}
 \langle P\rangle_\B(z) =
    \frac1{\left(2\pi\sqrt{-1}\right)^r}\int_{Z_\delta(\bl,q)\cap U(z)}
  \frac{P(\alpha_1,\dots,\alpha_n)\prod_{j=1}^r
    p_j\;\,d\mu_\Gamma^{\ga}}{(1-\kappa)
    \prod_{i=1}^n\alpha_i\prod_{j=1}^r(p_j-q_j)}. 
\end{equation}
\end{prop}

\begin{proof}
Consider the function \[
R=\frac{P(\alpha_1,\dots,\alpha_n)}
{D^\B(\alpha_1,\dots,\alpha_n)(1-\kappa)} \] and the differential
form
\[ \omega = \frac{dp_1}{p_1-q_1}\wedge\dots\wedge
\frac{dp_r}{p_r-q_r}\]
on $\comp$. Because of our assumptions, $\omega$ has a simple pole with
residue equal 1 at each of the points of the finite set $O_\B(z)$, and
the function $R$ is regular at these points. Our computation of the
Jacobian of the map $p$ shows that the divisors $\{u;\;p_j(u)=q_j\}$,
$j=1,\dots,r$ intersect transversally at these points, and thus for
small $\delta$ the set  $Z_\delta(\bl,q)$ consists of tiny tori, one
for each point of $O_\B(z)$ plus, possibly, some additional components
which we eliminate using the neighborhood $U(z)$ of $O_\B(z)$.

Then according to the usual integral representation of residues, we
have
\[ \langle P\rangle_\B(z) =
    \frac1{\left(2\pi\sqrt{-1}\right)^r}\int_{Z_\delta(\bl,q)\cap U(z)}
  \frac{P(\alpha_1,\dots,\alpha_n)\bigwedge_{j=1}^rdp_j}{(1-\kappa)
   D^\B(\alpha_1,\dots,\alpha_n)\prod_{j=1}^r(p_j-q_j)}.
 \]
Substituting the expressions from (2) and (3) into this formula, we
obtain (\ref{pb_integral}).
\end{proof}
\begin{remark}
1.  Note that, ``miraculously'', the differential form under the integral
sign here coincides with that in Proposition \ref{aside}.\\
2. We will show later that for  $z$ in a certain domain, the
conditions of the Proposition hold, moreover, $Z_\delta(\bl,q)$ is a
genuine cycle, i.e. it is localized in a small neighborhood of
$O_\B(z)$ and has no non-compact components. This last statement is
equivalent to the {\em properness} of the map $p$ defined in
\eqref{defp}, which, as we will see, turns out to be a subtle question.
\end{remark}

\section{Tropical calculations}
\label{sec:tropical}

In this section we only assume that $\A$ is a projective sequence in
$\Gamma_{\ga}^*$.  Recall that for each $\lambda\in \Gamma_\ga$ we
defined a rational function $p_\lambda(u)=\prod_{i=1}^n \alpha_i(u)^{
  \langle \alpha_i,\lambda\rangle }$ on $\ga_\C$, which is regular on
$U(\A)$; these functions $p_\lambda$ satisfy the relation
\begin{equation}
  \label{relation}
p_{\lambda_1+\lambda_2}(u)=p_{\lambda_1}(u) p_{\lambda_2}(u).
\end{equation}

\begin{definition}
Let $\xi\in \ga^*$. Define the set $$\ZZ(\xi)=\{u\in U(\A);\,\,
|p_\lambda(u)|=e^{-\langle  \xi, \lambda\rangle   } \text{ for all
}\lambda\in \Gamma_\ga\}.$$
\end{definition}

Fix a $\Z$-basis $(\lambda_1,\lambda_2,\ldots,\lambda_r)$ of the
lattice $\Gamma_\ga$, which is a positively oriented relative to our
chosen orientation of $\ga$, and introduce the notation $p_j$ for
$p_{\lambda_j}$. Then it follows from \eqref{relation} that the set
$\ZZ(\xi)$ is the subset of $U(\A)$ defined by the $r$ analytic
equations:
\begin{equation}
  \label{ranalytic}
|p_j(u)|=e^{- \langle \xi, \lambda_j\rangle }, j=1,\dots,r.
\end{equation}
In particular, $\ZZ(\xi)$ is an $r$-dimensional analytic subset of
$U(\A)$. We can orient $\ZZ(\xi)$ by the form $d\arg
p_1\wedge\dots\wedge d\arg p_r$. It is easy to see that this
orientation does not depend on the chosen positively oriented
basis.

The aim of this section is to prove that under some mild conditions
the cycle $\ZZ(\xi)$ is smooth, and compute its homology class in
$\comp$. 

Recall from \S\ref{sec:intersec}, that to each flag $F$ in
$\flaga$, one can associate a homology class $h(F)\in H_r(U(\A),\Z)$
(Definition \ref{def_hf}), a sign $\nu(F)$ given in Definition
\ref{goodflags}, and a sequence of  vectors
\[
\kappa^F_j=\sum\{\alpha_i;\;i=1,\dots,n,\, \alpha_i\in F_j\},
\]
given in  (\ref{kappaj}), which one can collect into a
sequence of $r$ vectors denoted by $\vkf$. By convention, we set
$\kappa^F_0=0$. Then for a flag $F\in \flaga$, introduce the non-acute
cone $\s(F,\A)$ generated by the non-negative linear combinations of
the elements $\{\kappa^F_j, j=1,\ldots,r-1\}$ and the line $\R
\kappa$:
$$\s(F,\A) =\sum_{j=1}^{r-1}\R^{\geq 0}\kappa_j+\R \kappa.$$

\begin{definition}\label{FF}
Let $\xi\in \ga^*$. We denote by $\FF$ the set of flags $F\in
\flaga$ such that $\xi\in\s(F,\A)$.
\end{definition}

Finally, recall from \S\ref{sec:intersec} that if $\xi$ is
sum-regular, then every flag $F\in \FF$ is proper, and thus has
$\nu(F)\neq0$.  The aim of this section is to prove the following theorem.
\begin{theorem}\label{Zxi}
  Let $\A$ be a projective sequence and let $\xi$ be a $\tau$-regular
  element in $\ga^*$ with $\tau$ sufficiently large. Then the set
  $\ZZ(\xi)$ is a smooth compact $r$-dimensional submanifold in
  $U(\A)$.  When oriented by the form $d \arg p_1\wedge d \arg
  p_2\wedge\cdots \wedge d \arg p_r$, it defines a cycle in $\comp$
  whose homology class $[\ZZ(\xi)]$ is given by
$$[\ZZ(\xi)]=\sum\nu(F) h(F),\quad {F\in \FF}.$$
\end{theorem}

The proof of the Theorem will start in \S\ref{sec:compact} and will
end with Proposition \ref{properto}. The compactness is contained in
Corollary \ref{compact}, the smoothness in Corollary \ref{smoothness}
and the computation of the homology class follows from Lemma
\ref{Lprop} and Proposition \ref{properto}.

\subsection{The tropical equations}
Our main tool of study of the set $\ZZ(\xi)$ is considering the
logarithm of the equations \eqref{ranalytic}, which can be written
as
\[\prod_{i=1}^n |\alpha_i(u) |^{\langle \alpha_i,\lambda_j\rangle
}=e^{-\langle \xi,\lambda_j\rangle },\; j=1,\dots,r.\] This idea
is related to tropical geometry \cite{Viro,Sturm}. The logarithmic
equations take the form
\[\sum_{i=1}^n \log |  \alpha_i(u)
|\langle \alpha_i, \lambda_j\rangle =-\langle
\xi,\lambda_j\rangle,\; j=1,\dots,r.
\]
This, in turn, can be written as a single vector equation:
\begin{equation}
  \label{vectoreq}
-\sum_{i=1}^n\log |\alpha_i(u)   |\, \alpha_i=\xi.
\end{equation}

Define the map $L:U(\A)\to \gd^*$ by
$$L(u)=-\sum_{i=1}^n\log | \alpha_i(u) | \omega^i.$$
Thus $L$ is a map from a real $2r$-dimensional space to an
$n$-dimensional one. If $u$ tends to $0$, then the point $L(u)$ tends
to $\infty$.

Recall from \S\ref{sec:prelim} that we denoted by $\mu$ the linear map
$\mu:\gd^*\to \ga^*$, which sends $\omega^i$ to $\alpha_i$. Then we
clearly have
 \begin{equation}
   \label{pl}
\left|p_j(u)\right|=e^{-\lr{\mu(L(u))}{\lambda_j}},
 \end{equation}
and thus  $u\in \ZZ(\xi)$ if and only if $\mu(L(u))=\xi$.
Another way to write this is that $\ZZ(\xi)=(\mu\circ L)^{-1}(\xi)$.

Our strategy is to separate the solution of \eqref{vectoreq} into
two parts, using that according to Lemma \ref{pl}, solutions to
\eqref{vectoreq} arise when the affine linear subspace
$\mu^{-1}(\xi)\subset\dgd$ of codimension $r$ intersects the image
$\iml$ of the map $L$.

Our next move is to give more precise information about the map $L$
and its image $\iml$.  Roughly, the idea is as follows. Since
\[
\log|\alpha_1(u)+\alpha_2(u)|\sim \max(\log|\alpha_1(u)|,\log
|\alpha_2(u)|)\] if $\alpha_1(u)$ and $\alpha_2(u)$ have different
orders of magnitude, we will be able to approximate the map $L$ with a
piecewise linear map from the $r$-dimensional space $U(\A)$ to
$\gd^*$.  Thus we will show that under some conditions the image
$\iml\subset\dgd$ is confined in a small neighborhood of a finite set
of affine linear subspaces of dimension $r$ which are transversal to
$\mu^{-1}(\xi)$.  This will allow us to describe the set $\ZZ(\xi)$
rather precisely.

Thus consider the affine subspace $\mu^{-1}(\xi)\subset \gd^*$. It
consists of the solutions $(t_1,t_2, \ldots, t_n)$ of the equation
\begin{equation}\label{system}
 \Eq(\xi):=\left\{\sum_{i=1}^n t_i\alpha_i= \xi\right\}.
\end{equation}
Often we will speak about the solutions of $\Eq(\xi)$ rather than
about the set $\mu^{-1}(\xi)$. Also, we will identify $\gd^*$ with
$\R^n$ whenever convenient, using the basis $\omega^i$.  Motivated by
Proposition \ref{tropical} below, we will be interested in a special
type of solutions of $\Eq(\xi)$ for which several of the coordinates
will be set equal.

\begin{definition} \label{flagdef}
  Let $F$ be a flag in $\flaga$ and $B=(B_1,B_2,\ldots,
  B_r)$ be a sequence of $r$ real numbers.  Define the point
  $\trop(F,B)=\sum_{i=1}^n t_i \omega^i\in \gd^*$ by the condition
  $t_i=B_j$ if $\alpha_i\in F_j\setminus F_{j-1}$. We will say that a
  solution $t\in \gd^*$ is an {\em  $F$-solution} of $\Eq(\xi)$
  \eqref{system} if $t$ is of the form $\trop(F,B)$ for some
     sequence $B$ of real numbers.
     \end{definition}

 Thus we see that $\trop(F,B)$ is a solution of $\Eq(\xi)$ if and only
if
\[ \sum_{j=1}^r
B_j(\kappa_j^F-\kappa_{j-1}^F)=\xi;\] this can be rewritten as :
\begin{equation}\label{Fsystem}
B_r\kappa+\sum_{j=1}^{r-1}(B_j-B_{j+1})\kappa^F_j=\xi.
 \end{equation}

 Recall that a flag is proper, if the elements $\kappa_j^F$ are
 linearly independent.  The following statement then clearly follows:
\begin{lemma}
  \label{indep}
For a proper flag $F$, there is exactly one  solution of
$\Eq(\xi)$ of the
  form $t(F,B)$. 
\end{lemma}
We denote this solution by $\sol(F,\xi)$.

For any flag $F$, consider  the system
\begin{equation}
  \label{Fsys}
\eqs(F) := \{t_b=t_c;\;\alpha_b,\alpha_c\in F_j\setminus
F_{j-1}\text{ for some }j\leq r\}.
\end{equation}

The solution set $G(F)$ of $\eqs(F)$ is a linear subspace of
$\dgd$ of dimension $r$ spanned by the vectors
\[s^{F,j}=\sum\{\omega^i;\;i=1,\dots,n, \alpha_i\in F_j\setminus
F_{j-1}\},\quad j=1,\dots,r.
\]
This subspace is transversal to the subspace $\mu^{-1}(0)$ of $\gd^*$
if and only $F$ is a proper flag. Indeed, the images of the vectors
$s^{F,j}$ under $\mu$ are equal to $\kappa_j^F-\kappa_{j-1}^F$, so
they span $\ga^*$ if and only if the vectors $\kappa_j^F$ do so.  Another
way to state Lemma \ref{indep} is to say that whenever $F$ is proper,
then $G(F)\cap \mu^{-1}(\xi)$ is non empty and consists of the single point
$\sol(F,\xi)$. Also note that if $\xi$ is regular with respect to
$\Sigma\A$ and there is an $F$-solution of $\Eq(\xi)$, then $F$ is
necessarily proper. Indeed the equation \eqref{Fsystem} implies that
$\xi$ is in the span of the vectors $\kappa_j^F$ belonging to
$\Sigma\A$.

Now we introduce a particular kind of $F$-solutions, which arise
naturally in our study of the set $\ZZ(\xi)$.

\begin{definition} \label{tropicaldef}
  Fix $\xi$, let $F$ be a flag in $\flaga$.  We will say that a
  solution $t\in \gd^*$ is a {\em tropical $F$-solution} of $\Eq(\xi)$
  \eqref{system} if $t$ is of the form $\trop(F,B)$ for some
    {\bf decreasing} sequence $B=(B_1,B_2,\ldots,
  B_r)$ of $r$ real numbers, that is with
  $B_1\geq\dots\geq B_r$. A solution of the equation $\Eq(\xi)$
   of the form $\trop(F,B)$ for some flag $F$ and some
  decreasing sequence $B$ will be called a {\em tropical solution} of
  $\Eq(\xi)$.  Finally, denote by $\FF$ the set of flags $F\in\flaga$
  for which $\Eq(\xi)$ has a tropical $F$-solution.
\end{definition}

 The following statement clearly follows from \eqref{Fsystem}.

\begin{lemma}\label{sfa}
  The flag $F$ belongs to $\FF$ if and only if $\xi\in\s(F,\A)$.
\end{lemma}

From now on, we will always assume that $\xi$ is regular with respect
to $\Sigma\A$. In particular, this implies that all flags $F\in \FF$
are proper.

\subsection{Compactness}\label{sec:compact}
Now we are ready to start the {\em Proof of Theorem \ref{Zxi}}.
We will show that when $\tau$ is sufficiently
large, then the cycle $\ZZ(\xi)$ is a disjoint union of compact smooth
components $\ZZF(\xi)$ associated to flags $F\in\FF$. The component
$\ZZF(\xi)$ will lie in an open set $U(F,N(\tau))$, where $N(\tau)$ is
increasing exponentially with $\tau$.  The sets $U(F,N)$ were defined
before Lemma \ref{ufn}.  The homology class of $\ZZF(\xi)$ will be a
generator of the $r$th homology of this set.

The first idea is that if $\xi$ is $\tau$-regular, then for every
$u\in \ZZ(\xi)$, there exists a flag $F\in\FF$ such that $L(u)$ is
close to the tropical solution $\sol(F,\xi)$. In fact, there is a
better approximation, as we show below.

In order to obtain a more precise estimate, we need to modify the
system $\eqs(F)$. As the space $F_j/F_{j-1}$ is 1-dimensional, for
two vectors $\alpha_b,\alpha_c\in F_j\setminus F_{j-1}$ there
exists a unique nonzero rational number $m_{bc}$ such that
$\alpha_b-m_{bc}\alpha_c\in F_{j-1}$.   Then let
\begin{equation}
  \label{modFsys}
\tS(F) := \{t_c-t_b=\log|m_{bc}|;\;\alpha_b,\alpha_c\in
F_j\setminus F_{j-1}\text{ for some }j\leq r\}.
\end{equation}
The solution set $\widetilde{G}(F)$ of $\tS(F)$ in $\dgd$ is an affine
$r$-dimensional subspace parallel to the solution set of $\eqs(F)$.
To be more specific, we construct concrete solutions of $\tS(F)$ as
follows.  Consider the ordered basis
$\bgam^F=(\gamma^F_1,\dots,\gamma^F_r)$ of $\dga$ that we introduced
earlier. It is such that $\{\gamma^F_m\}_{m=1}^j$ is a basis of $F_j$
for $j=1,\dots,r$. For $\alpha_b \in F_j\setminus F_{j-1}$, define the
rational number $m_b$ such that $\alpha_b -m_b\gamma_j\in F_{j-1}$.
Then the point $\sum_{i=1}^n-\log |m_i| \,\omega^i$ belongs to
$\tS(F)$, and $\tS(F)$ is the affine space parallel to $G(F)$ through
this point. This implies
\begin{lemma}
  \label{uniquesol}
  Let $F$ be a proper flag.  Then the solution spaces of the systems
  of linear equations $\Eq(\xi)$ and $\tS(F)$, $\mu^{-1}(\xi)$ and
  $\widetilde{G}(F)$ respectively, are transversal and of
  complementary dimensions.
\end{lemma}
It follows from this lemma that there is a unique common solution
of these equations; we denote this solution by $\ts(\xi,F)$.

The following is a key technical statement of this paper, which
justifies the validity of our tropical approximation.

Given a vector $t=(t_1,\dots,t_n)\in\dgd$, let
$\|t\|=\max_{i=1}^n|t_i|$ be the maximum norm of $t$.
\begin{prop}\label{tropical}
  There exist positive constants $\tau_0$, $c_0$ and $c_1$,
  which depend on $\A$ only, such that if $\tau\geq\tau_0$ and $\xi$
  is $\tau$-regular, then for every $u\in \ZZ(\xi)$ there exists a
  flag $F\in \FF$ such that
   $$\|L(u)-\ts(F,\xi)\|\leq c_0 e^{-c_1\tau}.$$
\end{prop}

We start the proof with two Lemmas.
\begin{lemma}\label{plus}
  Let $\bgam=(\gamma_1,\ldots,\gamma_r)$ be an ordered basis of a real
  vector space $V^*$ of dimension $r$. Let $u\in V_\C$ be such that
  $\gamma_j(u)\neq0$ for $j=1,\dots,r$, and set $B_j=-\log |
  \gamma_j(u)|$. There are positive constants $\lambda_0, \,c_0$ such
  that if $\lambda\geq \lambda_0$ and $B_j-B_{j+1}>\lambda$ for
  $j=1,\dots,r-1$, then the following holds: let $1\leq a\leq r$ and
  $\alpha=m_1\gamma_1+\cdots +m_a\gamma_a\in V^*$ a vector with
  $m_a\neq0$. Then $\alpha(u) \neq 0$, and furthermore,
  \begin{equation}
    \label{loga}
\left|\log |\alpha(u)/m_a|-\log |\gamma_a(u)|\right|
\leq c_0 e^{-\lambda}.
  \end{equation}
The constants $\lambda_0,c_0$ depend only on the data $(\bgam,
\alpha)$ and not of the element $u$ satisfying the hypothesis of the
Lemma.
\end{lemma}
\begin{remark}\label{Remark} Consider $\lambda$ a positive
  constant. Introduce the set
  $$U[\lambda]=\{u\in V_\C\,\,;\,\, |\log | \gamma_j(u)|- \log
  |\gamma_{k}(u) ||\geq \lambda, \text{ for all }j\neq k\}.$$ The lemma
  above implies that, provided $\lambda$ is sufficiently large, the
  distance between $\log | \alpha(u) |$ and the finite set
  $\{\log | \gamma_j(u)| , j=1,\ldots, r\}$ remains
  bounded as $u$ varies in the open set $U[\lambda]$.
\end{remark}
{\em Proof of Lemma \ref{plus}}:
   For
  $u\in V_\C$, with $\gamma_a(u) \neq 0$ we can write
$$   \alpha(u)    =m_a    \gamma_a(u)
\left(1+\sum_{k=1}^{a-1} \frac{m_k}{m_a}\frac{ \gamma_k(u) }{
\gamma_a(u)}\right).$$ This gives
$$\left(\frac{\alpha(u)}{m_a\gamma_a(u)}-1\right)=\sum_{k=1}^{a-1}
  \frac{m_k}{m_a}\frac{\gamma_k(u)}{ \gamma_a(u) }.$$

  Now assume $B_j- B_{j+1}\geq \lambda$ for $j=1,2,\ldots, r-1$.  Then
  for $k<a$, we have $|{ \gamma_k(u) }/{ \gamma_a(u)
  }|=e^{-(B_k-B_a)}\leq e^{-\lambda}$.

Define $\delta={\sum_{k=1}^{a-1} |m_k}/{m_a|}$. We obtain
$|\sum_{k=1}^{a-1} {m_k}{   \gamma_k(u) }/{ m_a\gamma_a(u)
}|\leq \delta e^{-\lambda}$ and thus
$$1-\delta e^{-\lambda}\leq \left|\frac{     \alpha(u)
}{m_a  \gamma_a(u) }\right |\leq 1+\delta e^{-\lambda}.$$

Let $\lambda_0$ be such that $\delta e^{-\lambda_0}\leq
\frac{1}{2}$, and assume that $\lambda\geq \lambda_0$. Then $
\alpha(u)\neq0 $, and taking logarithms and using the inequalities
$\log (1+x)\leq x,\;\log (1-x)\geq -2x $ for $0\leq x\leq
\frac{1}{2}$, we obtain

$$-2 \delta e^{-\lambda}\leq \log |{\alpha(u)}/{m_a}|-\log |\gamma_a(u)|
\leq \delta e^{-\lambda}.$$ Thus the estimate of the lemma holds
with $c_0=2\delta$, $\lambda_0=\log (2\delta)$.\qed
\smallskip

For $\alpha\in\dga$ and a basis $\brho\subset\Sigma\A$ of $\dga$,
we can write
\[ \alpha = \sum_{\gamma\in\brho}u_\gamma^\brho(\alpha)\gamma.\]
Introduce the constant
\[M(\A) = \max\{|u^\brho_\gamma(\alpha_i)|;\;i=1,\dots,n,\,\brho\text{
  basis of }\dga,\,
\brho\subset\Sigma\A,\,\gamma\in\brho\}.
\]
\begin{lemma}\label{distance}
  Let $\xi\in \a^*$ be a $\tau$-regular vector and $t=(t_1,t_2,\ldots,
  t_n)$ a solution of $\Eq(\xi)$. Then there exists an $r$-element
  subset $\sigma\subset\{1,\dots,n\}$, such that
  $$|t_i-t_j|\geq c_1\tau\text{ for } i,j\in \sigma,\, i\neq j,\text{
    where }c_1=\frac{1}{n M(\A)}.$$
\end{lemma}

\begin{proof}
  Let $\sigma\subset \{1,2,\ldots,n\}$ be a maximal subset satisfying the
  condition $|t_i-t_j|\geq c_1\tau$ for all $i,j\in \sigma, \,i\neq j$.
  We will arrive at a contradiction, assuming that the cardinality of
  $\sigma$ is strictly less than $r$.  For every $k\notin \sigma$, there
  exists $a(k)\in \sigma$ such that $|t_k-t_{a(k)}|<c_1\tau$. We can
  write

  \begin{multline*}
\xi=\sum_{i\in \sigma}t_i \alpha_i+\sum_{k\notin \sigma} t_k\alpha_k=
  \sum_{i\in \sigma} t_i \alpha_i+ \sum_{k\notin \sigma}
  (t_{k}-t_{a(k)})\alpha_k+t_{a(k)}\alpha_k\\
=\sum_{i\in \sigma} \sum_{k\notin \sigma, a(k)=i}t_i \left(\alpha_i+
\sum_{k\notin \sigma, a(k)=i}\alpha_k\right) +\sum_{k\notin \sigma}
  (t_{k}-t_{a(k)})\alpha_k.
      \end{multline*}

 Consider the set  $\brho_\sigma=\{\alpha_i+\sum_{k\notin \sigma,
    a(k)=i}\alpha_k;\;i\in\sigma\}$. This set is a subset of $\Sigma\A$,
  and it spans a vector space of dimension strictly less than $r$.
 By passing to a subset if necessary, we can assume that the set
 $\brho_\sigma\subset\Sigma\A$ is linearly independent. Let
 $\brho\subset\Sigma\A$ be a basis of $\dga$ containing $\brho_\sigma$.
 For an element $\gamma\in\brho\setminus\brho_\sigma$
  we can write the $\gamma$-coordinate
 of $\xi$ as
 \[ u^\brho_\gamma(\xi)=
 \sum_{k\notin \sigma} (t_k-t_{a(k)}) u^\brho_\gamma(\alpha_k).
\] Each number
$|u^\brho_\gamma(\alpha_k)|$ is less or equal than $M(\A)$, and
there are at most $n$ terms of this kind.  Thus
$|u_\gamma^\brho(\xi)|<\tau$. But this contradicts the
$\tau$-regularity of $\xi$.
\end{proof}

{\em Proof of Proposition \ref{tropical}}: Fix a $\tau$-regular vector
$\xi\in\a^*$, and let $u\in \ZZ(\xi)$. The vector $L(u)$ satisfies
$\Eq(\xi)$, and we can apply Lemma \ref{distance}; denote the index
subset guaranteed by the lemma by $\sigma(u)$.

Let $(B_1(u),\dots,B_r(u))$ be the set of numbers
$\{t_i;\;i\in\sigma(u)\}$ arranged in decreasing order, and if
$B_j(u)=t_i$, then we will write $\gamma_j$ for $\alpha_i$. Then we have
\begin{equation}
  \label{ufirst}
|\gamma_j(u)|=e^{-B_j(u)}\text{ with } B_j(u)-B_{j+1}(u)\geq
c_1\tau.
\end{equation}

Next assuming that $\tau$ is sufficiently large, we need to show that
the vectors $\{\gamma_j;\;j=1,\dots r\}$ are linearly independent. We
may use Lemma \ref{plus}. Indeed, assume that there is a linear
relation between these vectors. Let $\gamma_j$ and $\gamma_k$ be the
two vectors with the largest indices that have non-vanishing
coefficients in this relation.  Then according to \eqref{loga}, we
would have $|B_j(u)-B_k(u)|<m$ for some constant $m$ that only depends
on $\A$. This clearly cannot happen if $\tau$ is sufficiently large.

Next, denote by $F(u)$ the flag in $\flaga$ given by the sequences of
subspaces $F_j=\sum_{k=1}^j \C \gamma_k$, $j=1,\dots,r$.  Now we
estimate $\|L(u)-\ts(F(u),\xi)\|$. According to Lemma \ref{plus}, for
$\alpha_i\in F_j\setminus F_{j-1}$ we have
\begin{equation}
  \label{usecond}
|\log|\alpha_i(u)|-\log|\gamma_j(u)|-\log|m_i||\leq
c_2e^{-c_1\tau},
  \end{equation}
where $m_i$ is a constant such that $\alpha_i -m_i \gamma_j\in
F_{j-1}$, and $c_1,c_2$ are positive constants depending on $\A$
only.

Note that this equation implies that $F$ is in  $\FF$, provided
$\tau$ is sufficiently large. Indeed, the equation $\mu(L(u))=\xi$
implies that
\[
\left\|\xi-\left(-\log|\gamma_r(u)|\kappa+\sum_{j=1}^{r-1}
(\log|\gamma_{j+1}(u)|-\log|\gamma_{j}(u)|)\kappa_j^F\right)\right\|
\]
is uniformly bounded. This implies that for $\tau$ large the vectors
$\kappa_j^F$ have to be linearly independent, as otherwise $\xi$ would
not be $\tau$-regular. So $F$ is a proper flag; furthermore, for
$j=1,\ldots,r-1$, the number $\log |\gamma_{j+1}(u)|-\log
|\gamma_{j}(u)|=B_j(u)-B_{j+1}(u)$ is positive, bounded below by a
quantity that increases linearly with $\tau$. Thus $\xi$ is in the
cone $\s(F,\A)$ and  $F\in \FF$.

Let us look more closely at the equation \eqref{usecond}. The point
$\tilde t=\sum_{i=1}^n \tilde t_i \omega^i$ with coordinate $\tilde
t_i=-(\log |m_i|+ \log |\gamma_j(u)|)$ belongs to the affine space
$\tS(F)$ defined in \eqref{modFsys}.  The point $L(u)$ belongs to a
translate $G_u(F)$ of the linear subspace $G(F)\subset\dgd$, the
solution set of $\eqs(F)$. Thus this translate is at distance constant
times $e^{-c_1\tau}$ from the solution set $\widetilde{G}(F)$ of $\tS(F)$.
Since $L(u)$ also satisfies $\Eq(\xi)$, using Lemma \ref{uniquesol},
we see that $L(u)$, which is the intersection point of $G_u(F)$ and
the solution set of $\Eq(\xi)$, is not further from $\ts(F,\xi)$ than
a constant times $e^{-c_1\tau}$.  \qed

Using Proposition \ref{tropical} and its proof we can describe the
structure of $\ZZ(\xi)$ as follows.

Recall the definition and properties of the open sets $U(F,N)$,
and the constant $N_0$ given in Lemma \ref{ufn}. 
Introduce the sets
\[ \ZZF(\xi)=\ZZ(\xi)\cap U(F,N).\]

 \begin{cor}\label{compact} Let $N>N_0$.  Then there is a positive
 $\tau_0$ such that if $\xi$ is $\tau$-regular, with
 $\tau\geq \tau_0$, then $\ZZ(\xi)$ is the disjoint union of the compact
 sets $\ZZF(\xi)=\ZZ(\xi)\cap U(F,N)$ for $F\in\FF$.
\end{cor}
Note that in \eqref{ufirst} and \eqref{usecond}, at the cost of
changing the constants $c_1$ and $c_2$ we could replace the basis
vectors $\gamma_j$ from $\A$ by the basis $\bgam^F$ we chose for the
flag $F=F(u)\in\FF$. Then these equations, and the observation that
the map $L:\comp\ar\dgd$ is proper imply the statement. \qed

\subsection{Smoothness}
Next we turn to proving the smoothness of $\ZZF(\xi)$.  Fix a
$\Z$-basis $\bl=(\lambda_1,\lambda_2,\ldots,\lambda_r)$ of
$\Gamma_\ga$, and consider the map $p: U(\A)\to \C^{*r}$ given by
$$p(u)=(p_1(u),p_2(u),\ldots, p_r(u)),$$
where, as usual, $p_j$ stands for $p_{\lambda_j}$.

As  $\ZZ(\xi)$ is the inverse image of a smooth torus under that map
$p$, to prove that it is smooth, it is sufficient to show that the
Jacobian matrix of the map $p:\comp\ar \C^{*r}$ is
non-degenerate for $u\in U(F,N)$. According to Lemma \ref{DA} (1)
and (3), this reduces to the computation of
\begin{equation}
\label{td}
 \tD_{\A}(u) = \sum_{\sigma\in\bA}\vol_{\dga}(\bgam^\sigma)^2
  \prod_{i\in\sigma}\frac1{\alpha_i(u)}.
\end{equation}

Recall that in \S\ref{sec:intersec} we associated to each flag
$F\in\flaga$ a sequence of vectors $(\kappa_1^F,\dots,\kappa_r^F)$,
and we fixed a basis $\bgam^F$ satisfying
\[ d\gamma^F_1\wedge d\gamma^F_2\wedge \cdots \wedge d\gamma^F_r=
\,d\mu_\Gamma^{\ga}.\] We consider the $\gamma^F_j$s as
coordinates on $\ga_\C$ and to simplify our notation, we use
$u_j=\gamma^F_j(u)$ for $u\in\ga_\C$ and $j=1,\dots,r$.

\begin{prop} \label{df}
  Let $d(F)$ be the integer such that $\kappa_1^F\wedge \cdots \wedge
  \kappa_{r-1}^F \wedge \kappa=d(F)\,\,d\mu_\Gamma^{\ga}$.  Then for
  $N$ sufficiently large, we have
  $$\left|\tD_{\A}(u)\prod_{j=1}^ru_j-d(F)\right|\leq
  \frac{\const(\A)}{N}$$
  for any $u\in U(F,N)$.
\end{prop}
\begin{remark} Here and below we  use the same notation $\const(\A)$
  for several constants which depend only on $\A$.
\end{remark}

\begin{proof} Define $\bAF$ to be the set of those $\sigma\in\bA$ for which
 $\{\alpha_i, i\in\sigma\} \cap F_j$ has $j$ elements. Then the sum
 formula for $\tD_{\A}(u)$ is divided into two parts
 $\tD_{\A}(u)=\tD_\A^F(u)+R_\A^F(u)$, a dominant and a remainder
 term, where
$$\tD_\A^F(u)=\sum_{\sigma\in\bAF}\frac{\vol_{\A}(\bgam^\sigma)^2}
{\prod_{i\in\sigma}\alpha_i(u)  }$$ and

$$R_\A^F(u)=\sum \frac{\vol_{\A}(\bgam^\sigma)^2}
{\prod_{i\in \sigma}\alpha_i(u)},\quad \sigma\in\bA\setminus\bAF.
$$

Assuming $u\in U(F,N)$, we have the following basic estimates. For
$\alpha_i\in F_j\setminus F_{j-1}$, we have
\begin{equation}
  \label{basic}
\left|\frac{\alpha_i(u)}{u_j}-m_i\right|<\frac{\const(\A)}N,
\end{equation}
where $\alpha_i-m_i\gamma^F_j\in F_{j-1}$.
 This immediately leads to the estimate
\[ \left|R_\A^F(u)\prod_{j=1}^ru_j\right|<\frac{\const(\A)}{N}.
\]

Now we estimate  $\tD^F_\A(u)$. First, two simple linear algebraic
equalities:
\begin{lemma} \label{linalg}
\begin{enumerate}
\item $\prod_{i\in\sigma} m_i = \vol_{\A}(\bgam^\sigma)$ ,
\item $\sum_{\sigma\in\bAF}\vol_{\dga}(\bgam^\sigma) = d(F)$.
\end{enumerate}
\end{lemma}
\begin{proof}
The first equality follows from the fact that the matrix of
basis-change from $\{\alpha_i;\;i\in\sigma\}$ to the basis $\bgam^F$
is triangular, with the constants $m_i$ in the diagonal.  The second
one can be seen by expanding the sums
\[\kappa_j^F=\sum\{\alpha_i|\;\alpha_i\in F_j,\,
i=1,\dots,n\}\] in the exterior product $\kappa_1^F\wedge\dots
\wedge\kappa_r^F$. The non-vanishing terms will exactly correspond
to the sum on the left hand side of (2). 
\end{proof}

Now we can finish the proof of Proposition \ref{df}.
The sum defining $\tD_\A^F(u)$ is indexed by the elements
$\sigma\in\bAF$. The term corresponding to $\sigma$ multiplied by
$\prod_{j=1}^ru_j$ may be estimated as follows.
Using \eqref{basic} and the first equality in Lemma \ref{linalg}
we have
\[
\left|\prod_{j=1}^ru_j\;
\prod_{i\in\sigma}\frac1{\alpha_i(u)}-\frac1{\vol_{\dga}(\bgam^\sigma)}\right|<\frac{\const(\A)}{N}.
\]

Summing this inequality over $\sigma\in\bAF$ and simplifying the
fraction in each term, we obtain that
$$\left|\tD^F_{\A}(u)\prod_{j=1}^ru_j-\sum_{\sigma\in\bAF}
  \vol_{\dga}(\bgam^\sigma)\right|\leq \frac{\const(\A)}N.$$
Then applying
the second equality of Lemma \ref{linalg} completes the proof.
\end{proof}
As we observed earlier, (\ref{td}) is up to a nonzero multiple the
Jacobian of the map $p$. Then Corollary \ref{compact} together with
Proposition \ref{df} implies
\begin{cor}\label{smoothness}
For sufficiently large $\tau$ and $N$ the compact sets
$\ZZF(\xi)=\ZZ(\xi)\cap U(F,N)$,  $F\in\FF$, are smooth manifolds.
\end{cor}

\subsection{The homology class}
Now we turn to the computation of the homology class of the manifolds
 $\ZZF(\xi)$. 
 Introduce the torus
$$C(\xi)=\{(y_1,y_2,\ldots, y_r);\; |y_j|=e^{-\langle \xi, \lambda_j
  \rangle}\}\subset\C^{*r}.$$
If we orient $C(\xi)$ by the
differential form $d\arg y_1\wedge\dots\wedge d\arg y_r$, then its fundamental
class is a generator of $H_r(\C^{*r},\Z)$ over $\Z$.  Clearly, our set
$\ZZ(\xi)$ is the inverse image of $C(\xi)$ by the map
$p=(p_1,p_2,\ldots,p_r)$:
$$\ZZ(\xi)=\{u\in U(\A);\;
|p_1(u)|=e^{-\langle\xi,\lambda_1\rangle},\ldots,
|p_r(u)|=e^{-\langle\xi,\lambda_r\rangle}\}.$$

We can summarize what we have shown so far as follows.  Let $\tau$
be sufficiently large, positive and let $N=c_1e^{c_2\tau}$. Then
according to Corollary \ref{compact}, for a $\tau$-regular
vector $\xi$ the set $\ZZ(\xi)$ breaks up into finitely many
compact components $\ZZF(\xi)=\ZZ(\xi)\cap U(F, N)$, as $F$ varies
in $\FF$.  Since $\xi$ is $\tau$-regular, we have $d(F)\neq 0$ for
every flag in the family $\FF$.  Thus, according to \ref{df}, the
differential of the map $p$ does not vanish on $U(F,N)$ and thus
$\ZZF(\xi)$ is a smooth compact submanifold of $U(F,N)$,

What remains to prove Theorem \ref{Zxi} is that the homology class of
the oriented smooth manifold $\ZZF(\xi)$ in $H_r((U(F,N),\Z)$ is equal
to $\nu(F)h(F)$, where $h(F)$ is the fundamental class of the torus
$T_F(\ve)$ defined in \eqref{zefepsilon}. We will achieve this using a
deformation argument.

Recall that we have fixed an $F$-basis
$(\gamma_1^F,\gamma^F_2,\ldots, \gamma^F_r)$ of $\a^*$. Then for
$\alpha_i\in F_j\setminus F_{j-1}$, $i=1,\dots,n$, we can write
\begin{equation}
  \label{alphagamma}
\alpha_i= \sum_{k=0}^{j-1} m_{i,k} \gamma^F_{j-k},\;\text{ with }
m_i:=m_{i,0}\neq 0.
\end{equation}
Now we define a deformation  $\A_s^F$ of our sequence $\A$ as
follows:
\[\alpha^F_i(s,u)=\sum_{k=0}^{j-1} s^{k}m_{i,k}
u_{j-k}\quad\text{if }\alpha_i\in F_j\setminus F_{j-1},\]
where we again used the simplified notation $u_j=\gamma^F_j$.

In particular, we have $\alpha_i^F(1,\cdot)=\alpha_i$ and
$\alpha^F_i(0,\cdot)=m_{i} \gamma^F_j$. Using the estimate of Lemma
\ref{plus} we see that $\alpha^F_i(s,\cdot)$ does not vanish on
$U(F,N)$ for any $s\in [0,1]$ provided $N\geq N_0$.  This means that
we obtain a deformation of the map $p=(p_1,\dots,p_r)$ as well.
Define $p_j^F(s,u)=\prod_{i=1}^n
\alpha^F_i(s,u)^{\langle\alpha_i,\lambda_j\rangle}$.  Consider the map
$p^F: [0,1]\times U(F,N)\to \C^{*r},\,p^F(s,u)=(p^F_1(s,u),\ldots,
p^F_r(s,u))$, and let $\ZZF_s(\xi)= p^F(s,\cdot)^{-1}(C(\xi))$ be the
induced deformation of our cycle $\ZZ(\xi)$.

Similarly, we can define the map
$$L^F:[0,1]\times U(F,N)\to \dga,\quad L^F_s(u)=-\sum_{i=1}^n
\log|\alpha_i^F(s,u)|\,\omega^i.$$
Again, we have $\ZZF_s(\xi)=(\mu\circ L^F_s)^{-1}(\xi)$.

Then a  direct computation yields the following equalities:
 \begin{lemma} For $j=1,\dots,r$, we have $p^F_j(1,u)=p_j(u)$ and
\[p^F_j(0,u) = \prod_{i=1}^nm_i^{\lr{\alpha_i}{\lambda_j}} \prod_{l=1}^r
u_l^{\lr{\kappa^F_l-\kappa^F_{l-1}}{\lambda_j}}.
\]
Similarly, we have $L^F_1(u)=L(u)$ and
$$L^F_0(u)=-\sum_{j=1}^r \sum_{\alpha_i\in F_j\setminus F_{j-1}}
(\log|m_i|+\log|u_j|)\omega^i.$$
 \end{lemma}

Next, we compute the cycles $\ZZF_0(\xi)$.
\begin{lemma} \label{lmuzero}
  For a certain sequence of real numbers $\ve$, we have
\begin{equation}
  \label{Zzero}
\ZZF_0(\xi) = T_F(\ve)\subset U(F,N),
\end{equation}
where
\[ \ZZF_0(\xi) = p^F(0,\cdot)^{-1}(C(\xi))= (\mu\circ L^F_0)^{-1}(\xi),
\]
and the torus $T_F(\ve)$ was defined in \eqref{zefepsilon}.  In
addition, the orientation of the torus $\ZZF_0(\xi)$, induced by the
form $d\arg p^F_1(0,\cdot)\wedge\dots\wedge d\arg p^F_r(0,\cdot)$,
will coincide with the orientation of $T_F(\ve)$, induced by the form
$d\arg\gamma^F_1\wedge\dots\wedge d\arg\gamma^F_r$ exactly when $\nu(F)=1$.
\end{lemma}
\begin{proof}
  The fact that $\ZZF_0(\xi)$ is a torus immediately follows from the
  fact that each $p_j^F(0,\cdot)$ is a monomial in the linear forms
  $\gamma^F_j$, $j=1,\dots,r$. In fact, it is not hard to compute the
  sequence $\ve=(\epsilon_1,\dots,\epsilon_r)$: if $\alpha_i\in
  F_j\setminus F_{j-1}$, then $\epsilon_j=e^{-(t_i-\log|m_i|)}$, where
  $t_i$ is the $i$th component of $\ts(F,\xi)$.

To compare the orientations, observe that
\[ \frac{d\arg p^F_j}{d\arg
\gamma^F_l}=\lr{\kappa_l-\kappa_{l-1}}{\lambda_j}, \quad j=1,\dots, r.\]
This shows that the two orientations coincide exactly if the basis
$\vkf$ is oriented the same way as the basis $\bgam^F$. By definition
this happens exactly when $\nu(F)=1$.
\end{proof}

Thus we obtained a deformation of the cycle $\ZZF(\xi)$ to a cycle
which manifestly represents the homology class $\nu(F)h(F)\in
H_r(\comp,\Z)$.  To complete the proof of Theorem \ref{Zxi}, it
remains to show that the homology class of the cycles does not change
in this family.  This is fairly standard. The background for this
material is \cite{BT}.

If one has a proper smooth map between smooth manifolds $\pi:U\ar V$ with
$\dim U-\dim V = k$, then there is a natural grade-preserving
pull-back map
\[ \pi^*: H^\bullet_\cp(V,\Z)\longrightarrow H^\bullet_\cp(U,\Z)\]
on compactly supported cohomology. Via Poincar\'e duality this induces
a natural map
\[ \pi^*: H_\bullet(V,\Z)\longrightarrow H_{\bullet+k}(U,\Z).\]
on the homology groups. This has the property that when a compact
submanifold $S\subset V$ consists of regular values, then $\pi^*$
applied to the fundamental class of $S$ is exactly the fundamental
class of the manifold $\pi^{-1}(S)$.

These maps are homotopy invariant in the sense that if now $\pi$ is a
proper map from $[0,1]\times U$ to $V$, then the maps $\pi^*(0,\cdot)$
and $\pi^*(1,\cdot)$ are equal on the homology groups.

Our map $p:\comp\ar\C^{*r}$, and the closely related map $\mu\circ
L:\comp\ar \dga$, however, are {\em not} proper! Thus we need a slight
generalization of the pull-back maps. Since a map is proper if the
inverse image of compact sets is compact, we could say that the map
$\pi$ is {\em proper to} $V'$, for some open subset $V'\subset V$, if
the $\pi^{-1}(K)$ is compact for any compact $K\subset V'$.

\begin{lemma}
\label{Lprop}
Let $N>N_0$. Then for sufficiently large $\tau$, the map $\mu\circ L$
is proper to the set of $\tau$-regular elements of $\dga$. Moreover,
this map remains proper, when restricted to the set $U(F,N)$, for any
$F\in\flaga$. 
\end{lemma}
The Lemma follows from Proposition \ref{tropical} and Corollary
\ref{compact}.

Now we can reformulate in these homological terms what we are
computing. Consider the connected component of the set of
$\tau$-regular elements containing $\xi$.  Then we
are trying to show that the pull-back $(\mu\circ L)^*\eta_\xi$ of the
generator $\eta_\xi$ of the zeroth
homology of this component is equal to
$\sum_{F\in\FF}\nu(F)h(F)$. According to Corollary \ref{compact}, to
prove this, it
is sufficient to show that the pull-back of $\eta_\xi$ under the map
$\mu\circ L$ {\em restricted} to $U(F,N)$ is $\nu(F)h(F)$.
According to Lemma \ref{lmuzero}, the map 
\[ \mu\circ L^F(0,\cdot):U(F,N)\ar\dga\]
has this property, and we need to show the same for the map
\[ \mu\circ L^F(1,\cdot):U(F,N)\ar\dga.\]

Now this discussion explains that our deformation argument is
justified as long as we have the following.
\begin{prop}
  \label{properto}
  Let $N>N_0$, and $\xi\in\dga$ be a $\tau$-regular vector for $\tau$
  sufficiently large. Then the restricted map $\mu\circ
  L^F:[0,1]\times U(F,N)\ar\dga$ is
  proper to the set of $\tau$-regular elements.
\end{prop}
This statement is proved exactly the same way as we proved Corollary
\ref{compact}, the analogous statement for the map $\mu\circ L$.  It
is easy to see that the relevant constants are exactly the ones
appearing in the expressions of $\alpha_i^F(s,\cdot)$ via the basis
$\bgam^F$ in \eqref{alphagamma}. These constants are clearly uniformly
bounded as $s$ varies in $[0,1]$.  This completes the proof of the
proposition and that of Theorem \ref{Zxi} as well.  \qed

\section{The proof of the main results}
\label{sec:completion}

\subsection{The construction of the cycle for the JK-residue}
We start with an important observation.
\begin{prop}\label{near_zero}
  Let $\A$ be a projective sequence, and let $\gc$ be a chamber with
  $\kappa\in \overline \gc$. If $\xi\in \c$ be regular with respect to
  $\Sigma\A$, then all flags in $\FF$ are in fact in $\FFF$.
\end{prop}

\begin{proof} Let $F\in\FF$.  Equation $\Eq(\xi)$ given in
\eqref{Fsystem} reads as
\begin{equation}
  \label{xib}
\xi-B_r\kappa=\sum_{j=1}^{r-1} (B_j-B_{j+1})\kappa^F_{j}.
\end{equation}
The vector $\sum_{j=1}^{r-1} (B_j-B_{j+1})\kappa^F_{j}$ is a
positive linear combination of those elements of $\A$ which lie in the
subspace $F_{r-1}$. Thus the right hand side of \ref{xib} belongs to the closed set
$\con(\A)_\mathrm{sing}$. The vector $\xi$ is in $\gc$, and $\kappa$
is in $\overline{\gc}$, thus the ray $\xi+\R^{\geq0}\kappa$ does not touch
$\con(\A)_{sing}$. Therefore we must have $B_r>0$.
\end{proof}

Now we can formulate one of the central results of this paper: an
explicit construction of a real algebraic cycle which represents the
JK-residue.
\begin{theorem}\label{C} 
  Let $\A$ be a projective sequence and $\gc\subset \ga^*$ be a
  chamber such that $\kappa\in \overline{\gc}$, and fix an arbitrarily
  small neighborhood $U_0$ of the origin in $\ga_\C$.  If $\tau$ is
  sufficiently large, then for any $\tau$-regular $\xi\in\gc$, the set
  $\ZZ(\xi)$ is a smooth compact cycle in $U_0\cap\comp$ whose homology class
  equals the class $h(\gc)\in H_r(\comp,\Z)$ of the Jeffrey-Kirwan residue.
 \end{theorem}
\begin{proof}
  Indeed, we computed the homology class of $\ZZ(\xi)$ in Theorem
  \ref{Zxi}, and Theorem \ref{A} combined with
  Proposition \ref{near_zero} implies that this class is exactly the
  homology class realizing the Jeffrey-Kirwan residue.
  
  To prove the other statement of the Theorem, let $\xi$ be a
  $\tau$-regular vector. Then following the argument in the proof of
  Proposition \ref{near_zero}, we see that we must have $B_r\geq \tau$
  in \eqref{xib}. If $u$ is in $\ZZ(\xi)$, then it follows that $L(u)$
  is close to the point $\ts(F,\xi)$. This means that if $\alpha_i\in
  F_j\setminus F_{j-1}$, then $|\log|\alpha_i(u)|-B_j|$ is less than a
  constant depending on $\A$ only. As $B_j\geq B_r\geq\tau$, this
  means that we have
\[ |\alpha_i(u)|\leq \const(\A)\,e^{-\tau},\,\text{ for
}i=1,\dots,n.\]
This inequality implies the second statement of the Theorem. 
\end{proof}

\begin{remark} Formally, this construction only gives a representative
  of $h(\gc)$ in the case $\kappa\in\bar\gc$. Note, however,that for
  {\em any} chamber $\gc$ of $\A$, there exists a sequence $\A'$
  consisting of repetitions of the elements of $\A$, such that
  $\kappa'\in\bar\gc$, where $\kappa'$ is the sum of the elements of
  $\A'$. Applying the Theorem to $\A'$ will produce a representative
  for $h(\gc)$.
\end{remark}

\subsection{The proof of the conjecture} \label{sec:theproof}
Now we are ready to complete
the proof of Theorem \ref{main1}.  We recall our setup and introduce
some new notation. We have a projective, spanning sequence $\A$ and a
chamber $\gc$ containing $\kappa=\sumn\alpha_i$ in its closure. We
also picked a $\gc$-positive basis $\bl$. For $z\in\C^{*n}$, we denote
$z^{\lambda_j}$ by $q_j$, and introduce the vector $q=(q_1,\dots
q_r)\in\C^{*r}$.  For two vectors $\xi_1,\xi_2\in\dga$, write
$\xi_1\overl\xi_2$ if $\lr{\xi_1}{\lambda_j}<\lr{\xi_2}{\lambda_j}$
for $j=1,\dots,r$.

Fix a small vector $\eta\in\dga$ with the property that
$\lr\eta{\lambda_j}>0$ for $j=1,\dots,r$, i.e. $0\overl\eta$. Now we
pick $\xi\in\dga$ such that every vector $\zeta$ satisfying $\xi-\eta\overl\zeta\overl\xi$ is $\tau$-regular with $\tau$ sufficiently
large. ``Sufficiently large'' here means large enough to satisfy the
conditions of the statements we use in the course of the proof.

For any subset of $S\subset\{1,\dots,r\}$, define the torus
\[
T_S(\xi,\eta)=\left\{q\in\C^{*r};\;|q_j|=\begin{cases}
    \exp(-\lr{\xi}{\lambda_j})\text{ if 
    }j\notin S,\\  \exp(-\lr{\xi-\eta}{\lambda_j})\text{ if
    }j\in S\end{cases}\right\}
\]
with its standard orientation, let $Z_S(\xi,\eta)=p^{-1}T_S(\xi,\eta)$
be the inverse image of this cycle under the map $p$ (see
\eqref{defp}), and introduce the ring-like domain
\[ R(\xi,\eta) = \{q\in\C^{*r};\;
\lr{\xi-\eta}{\lambda_j}<-\log|q_j|<\lr\xi{\lambda_j},\;j=1,\dots,r\},\] on
whose edges this tori are located. Also, denote the associated domain
in $z$-space by $W(\xi,\eta)$:
\begin{equation}\label{condonz}
W(\xi,\eta) = \{z\in\C^{*n};\;\xi-\eta\overl \mu\circ L(z)\overl\xi\}.
\end{equation}
Thus $W(\xi,\eta)$ is the pull-back of the domain $R(\xi,\eta)$ under
the mapping which associates $q\in\C^{*r}$ to $z\in\C^{*n}$.

We will omit $(\xi,\eta)$ from the notation is this does not cause
confusion. For example we will use $Z_S$ instead of $Z_S(\xi,\eta)$.

Let the open set $U_0$ in Theorem \ref{C} be the set
$\{u\in\ga_\C;\;|\kappa(u)|<1\}$, and assume that $\tau$ is large
enough to satisfy the conditions of Theorem \ref{C} with this choice of
$U_0$. Now let    
\[U(\bl,q)=\{u\in\ga_\C;\;p^+_j(u)\neq q_jp^-_j(u)\}\subset\ga_\C,
\]
and recall our meromorphic $r$-form
\[ \Lambda =
\frac{P(\alpha_1,\dots,\alpha_n)\,\,d\mu_\Gamma^{\ga}} {(1-\kappa)
  \prod_{i=1}^n\alpha_i\prod_{j=1}^r (1-q_j/p_j)}
\]
on $\ga_\C$, which is regular on $U_0\cap\comp\cap U(\bl,q)$, and
depends on $z$.

Now our final argument may be broken up into the following 4
statements.
\begin{prop}\label{final} Let $z\in W(\xi,\eta)$ and assume that
  $\tau$ is sufficiently large. Then
\begin{enumerate}
\item $\int_{Z_{\emptyset}}\Lambda = \langle P\rangle_{\A,\gc}(z)$
\item $\int_{Z_S}\Lambda =0$ if $S\neq\emptyset$.
\item The cycle $\sum_S (-1)^{|S|} Z_S$, where $|S|$ denotes the
  number of elements of $S$,  is homologous in $U_0\cap\comp\cap
  U(\bl,q)$ to the cycle
 \[Z_\delta(q)= p^{-1}\{y=(y_1,\dots,y_r);\;|y_j-q_j|=\delta,\,
 j=1,\dots,r\},\] oriented by the form $d\arg(y_1-
 q_1)\wedge\dots\wedge d\arg(y_r- q_r)$.
\item $\int_{Z_\delta(q)}\Lambda = \langle P\rangle_{\B}(z).$
\end{enumerate}
\end{prop}
\begin{proof}
  We chose $\tau$ large enough in order to be able to apply Theorem
  \ref{C} to each of the cycles $Z_S$, $S\subset\{1,\dots,n\}$. Thus
  we know that the homology class of $Z_S$ in $\comp$ is $h(\gc)$ and
  that $Z_S\subset U_0$ for every $S\subset\{1,\dots,n\}$. The
  inequalities \eqref{condonz} defining $W(\xi,\eta)$ imply that for
  $u\in Z_\emptyset$ we have $|q_j|<|p_j(u)|$ for $j=1,\dots,r$, and
  thus we can apply Proposition \ref{aside}. This proves the first
  statement of the Proposition.

For $m=1,\dots,r$, denote by $Z_m$ the cycle $Z_{\{1,\dots,m\}}$. Since
we can permute the elements of the basis $\bl$, it is sufficient to
prove Statement (2) for these cycles. Now reversing the logic of the
proof of Proposition \ref{aside}, we can expand the differential form
$\Lambda$, taking into account that for $u\in Z_m$ and $z\in W(\xi,\eta)$ we
have $|p_j(u)|<|q_j|$ for $j=1,\dots,m$, and $|p_j(u)|>|q_j|$ for
$j=m+1,\dots,r$. We obtain the convergent expansion
\begin{multline*}
\int_{Z_m}\Lambda =
 \frac{(-1)^m}{\left(2\pi\sqrt{-1}\right)^r}\sum\int_{Z_m}\prod_{j=1}^m\frac{p_j^{l_j+1}}{q_j^{l_j+1}}
\prod_{j=m+1}^r\frac{q_j^{l_j}}{p_j^{l_j}}\cdot
\frac{P(\alpha_1,\dots,\alpha_n)\;\,d\mu_\Gamma^{\ga}}
 {(1-\kappa) \prod_{i=1}^n\alpha_i},\\
 \text{where the sum runs over }l_j\in\Z^{\geq0},\,j=1,\dots, r.
\end{multline*}
Since $Z_m$ represents the JK-residue, the terms of this series are
again of the form $\langle P\rangle_{\lambda,\A,\gc}$ with
$\lambda=\sum_{j=1}^r l_j\lambda_j$, where now $l_j<0$ for $j=1,\dots,
m$, and $l_j\geq0$ for $j=m+1,\dots,r$. Such an expression vanishes,
however, according to Proposition \ref{notincone}. This completes the
proof of the 2nd statement.

Clearly, (3) will follow from a similar statement formulated for the
cycles $T_S$:
\[ \sum_{S\subset\{1,\dots,n\}} (-1)^{|S|}T_S\text{ is homologous to }
\{y\in\C^{*r};\; |y_j-q_j|=\delta, j=1,\dots,r\}
\]
in the open set $\{y\in\C^{*r}; \; y_j\neq q_j\text{ for }j=1,\dots r\}$.
This may be proved by the standard inclusion-exclusion argument and is
left to the reader.

Finally, note that the cycle $Z_\delta(q)$ coincides with the cycle
$Z_\delta(\bl,q)$ introduced before Proposition \ref{B}. Then the
4th statement  will follow from Proposition \ref{B} as soon as we
check the technical conditions that we assumed there. We have done all
the groundwork for this; we just need to collect the necessary
information here. 

First, note that in Corollary \ref{compact} we show that $O_\B(z)\subset
U(F,N)$ for some $F\in\flaga$, in Proposition~\ref{DA} we compute the
Jacobian of the map $p$, and in Proposition~\ref{df} we show that
this Jacobian does not vanish on $O_\B(z)$. Thus we can conclude that
the set $O_\B(z)= p^{-1}(q)$ is finite. As
$D^\B(\alpha_1(u),\dots,\alpha_r(u))$ coincides with this Jacobian up
to a nonzero multiple, we see that it will not vanish on $O_\B(z)$.

Next, it follows from Lemma \ref{properto}, the map $p$ is proper to the
domain $R(\xi,\eta)$, and this eliminates the need for intersecting
with the small neighborhood $U(z)$ of $O_\B(z)$. Finally, note that we
already assumed that $\kappa(u)\neq1$ for any $u$ such that $p(u)\in
R(\xi,\eta)$, thus $1-\kappa$ will not vanish on $O_\B(z)$.
\end{proof}

Proposition \ref{final} proves the equality $\langle
P\rangle_{\A,\gc}(z)=\langle P\rangle_{\B}(z)$ for all $z\in
W(\xi,\eta)$, starting from the localized sum definition ~(\ref{pb}) of
$\langle P\rangle_{\B}(z)$. If we use the fact that this localized sum
is a toric residue, and thus it is a rational function of $z$, then we
can conclude that the two sides of (\ref{eq:main}) coincide whenever
the series in the left hand side converges. In view of of Lemma
\ref{convergence}, this implies the full statement of Theorem
\ref{main1}.


\newpage

\begin{thebibliography}{99}



\bibitem{Bat1} V. V. Batyrev, {\em Dual polyhedra and mirror symmetry
    for Calabi-Yau hypersurfaces in toric varieties}, J. Alg. Geom.
    {\bf 3} (1994), 493-535.
  \bibitem{BM} V. V. Batyrev, E. N. Materov, {\em Toric residues and
      mirror symmetry}, Dedicated to Yuri I. Manin on the
      occasion of his 65th birthday, Mosc. Math. J. {\bf 2} (2002),
      no. 3, 435-475.
\bibitem{BT} R. Bott, L. Tu, {\em Differential Forms in Algebraic
    Topology}, Springer-Verlag, New York, 1982.
\bibitem{B-V} M. Brion, M. Vergne, {\em Arrangement of hyperplanes I~: Rational
functions and Jeffrey-Kirwan residue}, Ann. Scient. {\'E}c. Norm.
Sup., {\bf 32} (1999) 715-741.
\bibitem{CCD} E. Cattani, D. Cox, A. Dickenstein, {\em Residues in
    toric varieties}, Compositio Math. {\bf 108} (1997), 35-76.
\bibitem{CDS} E. Cattani, A. Dickenstein, B. Sturmfels. {\em Residues
    and resultants}, J. Math Sci. Univ. Tokyo {\bf 5} (1998), no.1,
    119-148.
\bibitem{Cox} D. Cox, {\em Toric residues}, Arkiv f\"or Matematik {\bf
    34} (1996), 73-96.
\bibitem{Danilov} V. I. Danilov, {\em The geometry of toric
    varieties}, Russian Math. Surveys {\bf 33} (1978), 95-154.
\bibitem{Fulton} W. Fulton, {\em Introduction to Toric Varieties},
  Princeton University Press, Princeton, 1993.
\bibitem{GKZ} I. M. Gelfand, M. M. Kapranov, A. V. Zelevinsky, {\em
    Discriminants, Resultants and Multidimensional Determinants},
    Birkh\"auser, Boston, 1994.
\bibitem{GH} P. Griffiths, J.  Harris, {\em Principles of Algebraic
    geometry}, John Wiley, New York, 1978.
\bibitem{HS} T. Hausel, B, Sturmfels, {\em Toric hyperK\"ahler
    varieties}, Doc. Math {\bf 7} (2002), 495-534.
\bibitem{JK} L. Jeffrey, F. Kirwan, {\em Localization for nonabelian
    group actions}, Topology {\bf 34} (1995) 291-327.
\bibitem{MP} D. Morrison, R. Plesser, {\em Summing the instantons:
    Quantum cohomology and mirror symmetry in toric varieties},
    Nucl. Phys. B {\bf 440} (1995), 279-354.
\bibitem{OT} P. Orlik, H. Terao, {\em Arrangements of Hyperplanes},
    Grundlehren der Mathematischen Wissenschaften. Vol 300,
    Springer-Verlag.  Berlin, 1992.
\bibitem{asz_ir} A. Szenes, {\em Iterated residues and Bernoulli
      Polynomials}, Internat. Math. Res. Notices {\bf 18} (1998)
    937-956.
\bibitem{VS} V. V. Schechtman, A. N. Varchenko, {\em Arrangements of
      hyperplanes and Lie algebra homology}, Invent. Math., {\bf 106}
    (1991), 139-194.
\bibitem{Sturm} B. Sturmfels, Solving Systems of Polynomial Equations,
  CBMS Regional Conference Series in Math., no.97, AMS, Providence, RI, 2002.
\bibitem{Viro} O. Viro, {\em Dequantization of real algebraic geometry on a
  logarithmic paper}, Proc. 3rd European Congress of Mathematics 2000,
  vol. I, Progr. Math. {\bf 201}, 135-146, Birkh\"auser, Basel, 2001. 
\bibitem{W} E. Witten, {\em Two-dimensional gauge theory revisited},
  J. Geom. Phys. {\bf 9} (1992), no. 4, 303-368.




\end{thebibliography}
\end{document}